\documentclass[12pt]{extarticle}
\usepackage[a4paper, total={6.5in, 10in}]{geometry}
\usepackage{amsmath}
\usepackage{amsfonts}
\usepackage{amssymb}
\usepackage{amsthm}
\usepackage{mathrsfs}
\usepackage{graphicx}
\usepackage{float}
\usepackage[margin=1in]{caption}
\usepackage{enumerate}
\usepackage{enumitem}
\usepackage{setspace}
\usepackage{tikz}
\usepackage{tikz-cd}
\usepackage{setspace}
\usepackage{changepage}

\usepackage{amsmath,amscd}
\RequirePackage{amsmath}
\RequirePackage{amsfonts}
\RequirePackage{amssymb}
\RequirePackage{mathrsfs}
\RequirePackage{float}
\RequirePackage{tikz}
\RequirePackage{tikz-cd}

\newcommand{\R}{\mathbb{R}}
\newcommand{\F}{\mathbb{F}}
\newcommand{\Z}{\mathbb{Z}}
\newcommand{\Q}{\mathbb{Q}}
\newcommand{\C}{\mathbb{C}}

\renewcommand{\.}{\dots}
\newcommand{\+}{+\dotsm+}

\newcommand{\injec}{\hookrightarrow}

\makeatletter
\newcommand*{\da@rightarrow}{\mathchar"0\hexnumber@\symAMSa 4B }
\newcommand*{\da@leftarrow}{\mathchar"0\hexnumber@\symAMSa 4C }
\newcommand*{\xdashrightarrow}[2][]{%
	\mathrel{%
		\mathpalette{\da@xarrow{#1}{#2}{}\da@rightarrow{\,}{}}{}%
	}%
}
\newcommand{\xdashleftarrow}[2][]{%
	\mathrel{%
		\mathpalette{\da@xarrow{#1}{#2}\da@leftarrow{}{}{\,}}{}%
	}%
}
\newcommand*{\da@xarrow}[7]{%
	\sbox0{$\ifx#7\scriptstyle\scriptscriptstyle\else\scriptstyle\fi#5#1#6\m@th$}%
	\sbox2{$\ifx#7\scriptstyle\scriptscriptstyle\else\scriptstyle\fi#5#2#6\m@th$}%
	\sbox4{$#7\dabar@\m@th$}%
	\dimen@=\wd0 %
	\ifdim\wd2 >\dimen@
	\dimen@=\wd2 %   
	\fi
	\count@=2 %
	\def\da@bars{\dabar@\dabar@}%
	\@whiledim\count@\wd4<\dimen@\do{%
		\advance\count@\@ne
		\expandafter\def\expandafter\da@bars\expandafter{%
			\da@bars
			\dabar@ 
		}%
	}%  
	\mathrel{#3}%
	\mathrel{%   
		\mathop{\da@bars}\limits
		\ifx\\#1\\%
		\else
		_{\copy0}%
		\fi
		\ifx\\#2\\%
		\else
		^{\copy2}%
		\fi
	}%   
	\mathrel{#4}%
}
\makeatother

\usepackage{amsmath,amscd}
\RequirePackage{amsmath}
\RequirePackage{amsfonts}
\RequirePackage{amssymb}
\RequirePackage{mathrsfs}
\RequirePackage{float}
\RequirePackage{tikz}
\RequirePackage{tikz-cd}

\renewcommand{\F}{\mathscr{F}}

\newcommand{\T}{\textnormal{T}}

\newcommand{\eval}{\operatorname{eval}}			% random useful operator names

\renewcommand{\bar}{\overline}					% replaces \overline with \bar

\usepackage{geometry}
\usepackage{lipsum}
\usepackage{booktabs}
\usepackage{multirow}
\usepackage{siunitx}
\usepackage{wrapfig}
\usepackage{nicefrac}

\renewcommand{\bar}{\overline}

\newcommand{\sing}{\textnormal{sing}}
\newcommand{\reg}{\textnormal{reg}}

\renewcommand{\l}{\ell}

\renewcommand{\T}{\mathbb{T}}

\newcommand{\D}{\mathcal{D}}

\newcommand{\mg}{\textnormal{mg}}
\renewcommand{\F}{\mathcal{F}}

\newcommand{\supp}{\textnormal{supp}}

\begin{document}
	
	\begin{center}
		{\bf ON SPECTRA OF METRIC TREE GRAPHS}\\~\\ {Tyler Chamberlain}\\
	\end{center}
	
	\begin{adjustwidth}{30pt}{30pt}
		{\small Abstract.} The secular manifold $\Sigma_G$ and its singularities are intimately related to the spectra of metric graphs $(G,\ell)$. In this paper, we present a complete description of the singular locus for tree graphs, and confirm that it agrees with a conjecture of Colin de Verdi\`ere. We also discuss numerous applications toward studying the behavior of metric tree graphs.\\
	\end{adjustwidth}
	
	\noindent {\bf 1. Introduction.} Lenoid Friedlander in 2005 first proved that for a generic choice of $\ell$, the metric graph $(G,\ell)$ has a simple spectrum, so long as $G$ is not the circle [Fr]. His approach was purely analytic and made no explicit reference to the secular manifold $\Sigma_G$. Since, the geometry of $\Sigma_G$ has been studied much more extensively, for instance by Colin de Verdi\`ere, Kurasov, and Sarnak [CdV] [K-S]. And recently, Lior Alon strengthened Friedlander's result by employing algebraic methods [Al]. Of primary concern is its singular locus; any bound on its codimension in $\Sigma_G$ translates to a new and improved genericity statement. In fact, Alon's result uses only that $\Sigma_G$ is nonreduced for a non-circle graph $G$, so in other words, the singular locus has codimension one or more in $\Sigma_G$. It has been conjectured by Colin de Verdi\`ere that this bound can be improved to at least two when $\Sigma_G$ is irreducible [CdV].  
	
	The goal of this paper is to confirm Colin de Verdi\`ere's conjecture for tree graphs. In fact, we present an explicit description of $\Sigma^{(m)}_G$ in terms of the zero loci of other secular polynomials. Here, $\Sigma_G^{(m)}$ denotes the collection of points in the secular manifold with multiplicity $\geq m$. The key insight is the so-called {\it smoothness criteria} for trees (Lemma 3.1), which establishes conditions for which a point $z$ on $\Sigma_G$ would be smooth. Many of the major obsticles were purely notational, and as a result, we introduce a fair amount of machinery.
	
	In section 4, we consider some basic applications of this characterization for the study of generic eigenfunctions and uniformly discrete metric trees. We do prove the following surprising result:\\
	
	{\noindent \bf \footnotesize THEOREM.} {{\it Fix a tree graph $G$. Then for a generic choice of rationally dependent $\ell$, $(G, \ell)$ has simple spectrum.}}\\

	%Essential to studying generic eigenfunctions on metric graphs $(G,\ell)$ is the secular (or determinant) manifold $\Sigma_G$. Many non-generic properties of eigenfunctions code into systems of polynomial equations, with solutions lying strictly within $\Sigma_G$. 
	
	%Spectral geometry is traditionally done in the setting of Riemannian manifolds, with or without boundary. At times, however, we permit singularities, i.e. points where a local Euclidean structure cannot be extended in any meaningful manner. The collection of singular points is usually subject to either topological or dimensional constraints. In the one-dimensional case, such spaces are called {\it metric graphs}; the reason being that they appear (up to isometry) as a topological graph $\G$ with a collection $\l$ of metrics on its edges. The singularities are isolated and correspond to the vertices of $\G$ with degree larger than one. Like in the smooth case, we may inquire about the spectra of metric graphs, although for this to have any geometric relevance, concessions must be made at the vertices, or boundary, of $\G$. In particular, functions are required to satisfy either Neumann or Dirichlet conditions at such points. 
	
	\noindent {\bf 2. Preliminaries.} Every 1-d singular Riemannian manifold exists as a metric graph $(G,\l)$. The choice of graph $G$ is not unique; there is no topological distinction between a single edge and two connected at a vertex. The current literature nonetheless sports a very combinatorial approach to these spaces. Both the secular polynomial $P_G$ and manifold $\Sigma_G$ are combinatorial (not topological) invariants of the graph. A majority circumvent the issue by disallowing degree two vertices altogether, ensuring a unique choice of $G$. Boundary conditions at the vertices of $G$ follow a similar convention: the presence of a Dirichlet vertex acts to `disconnect' the graph at that point, and as such, the common practice is to permit Dirichlet conditions only at endpoints (degree 1 vertices), thereby keeping the underlying topology intact. However, when studying the singular locus, we must confront some rather unusual graphs, possibly with degree two vertices or Dirichlet conditions which do not discriminate\textemdash that is to say, they could occur at any vertex of $G$. The bottom line is that we cannot exclude these unpleasant constructions. Instead, we develop an apparatus to deal with them appropriately. This section serves to do so, in addition to providing a standard treatment of metric graphs.

	%The graph $G$, however, is not unique since there is no topological distinction between a single edge and two connected at one vertex. The current literature nonetheless sports a very combinatorial approach to these spaces\textemdash the secular polynomial $P_G$ and manifold $\Sigma_G$ are not solely determined by the topology of $G$. The majority circumvent the issue by disallowing degree two vertices altogether; this ensures the uniqueness of $G$. Boundary conditions at the vertices of $G$ follow a similar convention. The presence of a Dirichlet vertex acts to `disconnect' the graph at that point, and as such, the common practice is to permit Dirichlet conditions only at endpoints (degree 1 vertices), thereby keeping the underlying topology intact. However, in our quest to understand the secular manifold $\Sigma_G$, we are forced to confront some rather unusual graphs, possibly with degree two vertices or Dirichlet conditions which do not discriminate\textemdash that is to say, they could occur at any vertex of the graph. The bottom line is that we cannot exclude these unpleasant constructions, and instead, we must develop an apparatus to deal with them appropriately. This section serves to do so, in addition to providing a standard treatment of metric graphs.
	
	{\it I. The Graph as a Space.} Let $G$ be a finite graph with enumerated edges $e_j$. A metric $\ell$ on $G$ is a tuple of positive real numbers, assigning edge $e_j$ length $\ell_j$. The metric graph $(G, \ell)$ is endowed with a Riemannian structure. Singularities of the space occur at vertices of degree $>1$, where a singularity is said to be {\it removable} if it has degree exactly two. The boundary $\partial G$ is the vertex set of the graph. Regular points of the boundary are in correspondence with the endpoints of the graph. Each edge is parametrized by intervals of the form $[0,\ell_j]$. A cover of $G$ by these charts gives a natural orientation of the graph and vice versa. In general, we pick an arbitrary orientation. Following [Fr], we write $e_j\succ v$ (resp. $e_j \prec v$) to denote that a vertex $v$ is a target (resp. source) of edge $e_j$. 
	
		\begin{center}
		\tikzset{every picture/.style={line width=0.75pt}} %set default line width to 0.75pt        
		
		\begin{tikzpicture}[x=0.75pt,y=0.75pt,yscale=-1,xscale=1]
			%uncomment if require: \path (0,480); %set diagram left start at 0, and has height of 480
			
			%Straight Lines [id:da012092636666941381] 
			\draw [line width=0.75]    (162,161) -- (92,161) ;
			%Straight Lines [id:da9531717350819835] 
			\draw [line width=0.75]    (243,191) -- (173.35,191) ;
			\draw [shift={(171,191)}, rotate = 180] [color={rgb, 255:red, 0; green, 0; blue, 0 }  ][line width=0.75]      (0, 0) circle [x radius= 3.35, y radius= 3.35]   ;
			\draw [shift={(243,191)}, rotate = 180] [color={rgb, 255:red, 0; green, 0; blue, 0 }  ][fill={rgb, 255:red, 0; green, 0; blue, 0 }  ][line width=0.75]      (0, 0) circle [x radius= 3.35, y radius= 3.35]   ;
			%Straight Lines [id:da6074280296866537] 
			\draw [line width=0.75]    (170.32,188.75) -- (162.68,163.25) ;
			\draw [shift={(162,161)}, rotate = 253.3] [color={rgb, 255:red, 0; green, 0; blue, 0 }  ][line width=0.75]      (0, 0) circle [x radius= 3.35, y radius= 3.35]   ;
			\draw [shift={(171,191)}, rotate = 253.3] [color={rgb, 255:red, 0; green, 0; blue, 0 }  ][line width=0.75]      (0, 0) circle [x radius= 3.35, y radius= 3.35]   ;
			%Straight Lines [id:da11327477803255159] 
			\draw [line width=0.75]    (162,161) -- (153,131) ;
			\draw [shift={(153,131)}, rotate = 253.3] [color={rgb, 255:red, 0; green, 0; blue, 0 }  ][fill={rgb, 255:red, 0; green, 0; blue, 0 }  ][line width=0.75]      (0, 0) circle [x radius= 3.35, y radius= 3.35]   ;
			\draw [shift={(162,161)}, rotate = 253.3] [color={rgb, 255:red, 0; green, 0; blue, 0 }  ][fill={rgb, 255:red, 0; green, 0; blue, 0 }  ][line width=0.75]      (0, 0) circle [x radius= 3.35, y radius= 3.35]   ;
			%Straight Lines [id:da32167156128954444] 
			\draw [line width=0.75]    (180,221) -- (171.68,193.25) ;
			\draw [shift={(171,191)}, rotate = 253.3] [color={rgb, 255:red, 0; green, 0; blue, 0 }  ][line width=0.75]      (0, 0) circle [x radius= 3.35, y radius= 3.35]   ;
			\draw [shift={(180,221)}, rotate = 253.3] [color={rgb, 255:red, 0; green, 0; blue, 0 }  ][fill={rgb, 255:red, 0; green, 0; blue, 0 }  ][line width=0.75]      (0, 0) circle [x radius= 3.35, y radius= 3.35]   ;
			%Straight Lines [id:da9033011953275529] 
			\draw [line width=0.75]    (63,161) -- (92,161) ;
			\draw [shift={(63,161)}, rotate = 0] [color={rgb, 255:red, 0; green, 0; blue, 0 }  ][fill={rgb, 255:red, 0; green, 0; blue, 0 }  ][line width=0.75]      (0, 0) circle [x radius= 3.35, y radius= 3.35]   ;
			%Straight Lines [id:da9829742203111733] 
			\draw [line width=0.75]    (243,191) -- (271.65,191) ;
			\draw [shift={(274,191)}, rotate = 0] [color={rgb, 255:red, 0; green, 0; blue, 0 }  ][line width=0.75]      (0, 0) circle [x radius= 3.35, y radius= 3.35]   ;
			%Straight Lines [id:da5443595058792894] 
			\draw [line width=0.75]    (92,161) -- (83.68,133.25) ;
			\draw [shift={(83,131)}, rotate = 253.3] [color={rgb, 255:red, 0; green, 0; blue, 0 }  ][line width=0.75]      (0, 0) circle [x radius= 3.35, y radius= 3.35]   ;
			%Straight Lines [id:da22487501563302792] 
			\draw [line width=0.75]    (252,221) -- (243,191) ;
			\draw [shift={(243,191)}, rotate = 253.3] [color={rgb, 255:red, 0; green, 0; blue, 0 }  ][fill={rgb, 255:red, 0; green, 0; blue, 0 }  ][line width=0.75]      (0, 0) circle [x radius= 3.35, y radius= 3.35]   ;
			\draw [shift={(252,221)}, rotate = 253.3] [color={rgb, 255:red, 0; green, 0; blue, 0 }  ][fill={rgb, 255:red, 0; green, 0; blue, 0 }  ][line width=0.75]      (0, 0) circle [x radius= 3.35, y radius= 3.35]   ;
			%Curve Lines [id:da2409524489747361] 
			\draw [line width=0.75]    (92,161) .. controls (102.78,184.52) and (136.61,191.71) .. (169.02,191.05) ;
			\draw [shift={(171,191)}, rotate = 358.26] [color={rgb, 255:red, 0; green, 0; blue, 0 }  ][line width=0.75]      (0, 0) circle [x radius= 3.35, y radius= 3.35]   ;
			%Straight Lines [id:da11532997198394201] 
			\draw [line width=0.75]    (243,191) -- (234,161) ;
			%Straight Lines [id:da9997459367727586] 
			\draw [line width=0.75]    (234,161) -- (225.68,133.25) ;
			\draw [shift={(225,131)}, rotate = 253.3] [color={rgb, 255:red, 0; green, 0; blue, 0 }  ][line width=0.75]      (0, 0) circle [x radius= 3.35, y radius= 3.35]   ;
			\draw [shift={(234,161)}, rotate = 253.3] [color={rgb, 255:red, 0; green, 0; blue, 0 }  ][fill={rgb, 255:red, 0; green, 0; blue, 0 }  ][line width=0.75]      (0, 0) circle [x radius= 3.35, y radius= 3.35]   ;
			%Straight Lines [id:da11120615502916542] 
			\draw [line width=0.75]    (234,161) -- (265,161) ;
			\draw [shift={(265,161)}, rotate = 0] [color={rgb, 255:red, 0; green, 0; blue, 0 }  ][fill={rgb, 255:red, 0; green, 0; blue, 0 }  ][line width=0.75]      (0, 0) circle [x radius= 3.35, y radius= 3.35]   ;
			\draw [shift={(234,161)}, rotate = 0] [color={rgb, 255:red, 0; green, 0; blue, 0 }  ][fill={rgb, 255:red, 0; green, 0; blue, 0 }  ][line width=0.75]      (0, 0) circle [x radius= 3.35, y radius= 3.35]   ;
			%Straight Lines [id:da06262557196625163] 
			
			%Straight Lines [id:da20161292975614176] 
			%Straight Lines [id:da6482043713785761] 
			%Straight Lines [id:da4487694374085558] 
			%Straight Lines [id:da570524739607791] 
			%Straight Lines [id:da7132939524127442] 
			%Straight Lines [id:da5088190049152246] 
			%Straight Lines [id:da9327828103504485] 
			%Straight Lines [id:da646531038543102] 
			%Curve Lines [id:da055389892722882195] 
			%Curve Lines [id:da2617370517208757] 
			\draw    (172.43,189.06) .. controls (187.74,168.59) and (200.05,161.97) .. (234,161) ;
			\draw [shift={(234,161)}, rotate = 358.36] [color={rgb, 255:red, 0; green, 0; blue, 0 }  ][fill={rgb, 255:red, 0; green, 0; blue, 0 }  ][line width=0.75]      (0, 0) circle [x radius= 3.35, y radius= 3.35]   ;
			\draw [shift={(171,191)}, rotate = 306.03] [color={rgb, 255:red, 0; green, 0; blue, 0 }  ][line width=0.75]      (0, 0) circle [x radius= 3.35, y radius= 3.35]   ;\end{tikzpicture}
		\begin{adjustwidth}{75pt}{75pt}
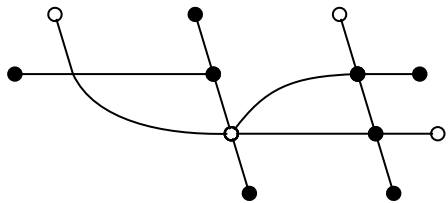
			{\small {\bf Figure 2.1.} A non-compact graph $G$ along with boundary $\partial G$. Closed (resp. open) points symbolize vertices (resp. not) contained in $G$.}
		\end{adjustwidth}
	\end{center}
	
	Boundary conditions are embedded in the underlying topological space. Constructing a space as in the paragraph above produces a graph with Neumann conditions everywhere. Dirichlet conditions at a vertex are specified by its removal from the space; for instance, the incomplete graph $G\setminus v$ refers the space $G$ with Dirichlet conditions at $v$. Conversely, we may start with an incomplete graph $G$ and specify its boundary $\partial G$ as the collection of vertices (see Figure 2.1). If omitted, $\partial G$ is assumed to contain all the vertices of $G$ in addition to a unique vertex for each incomplete edge. Moving forward, we will simply refer to these spaces as {\it graphs} to remain consistent, although they need not be compact spaces. The primary advantage of this approach is that boundary conditions are inherited by subspaces. A {\it subgraph} $H\subset G$ is a subspace of $\bar G$ with $\partial H\subset \partial G$. Thus, the closure of $H$ is a legitmate subgraph of the closure of $G$. 
	
	{\it II. Eigenspaces of a Metric Graph.} Functions on a metric graph $(G,\ell)$ are complex-valued functions on the closure $\bar G$; they are continuous everywhere and smooth at interior points. Each function $\phi$ appears in coordinates on edge $e_j$ as $\phi_j$. Likewise the Laplace operator $\Delta$ restricts to its familar 1-d expression $-d^2/dx_j^2$. The domain of $\Delta$ consists of functions vanishing at every vertex not in $G$ and or with the Neumann condition $$\sum_{e_j\succ v} \frac{d\phi_j}{dx_j}(v) = \sum_{e_j\prec v} \frac{d\phi_j}{dx_j}(v)$$ at vertices in $G$. An eigenfunction $\phi$ of $\Delta$ satisfies the differential equations $$\frac{d^2\phi_j}{dx_j^2} + k^2 \phi_j = 0$$ for all $j$ and some real $k > 0$. For the sake of brevity, such a function is called a {\it k-eigenfunction}. It appears in coordinates as $$\phi_j(x_i) = a_j \exp(ikx_j) + b_j z_j\exp(-ikx_j)$$ for some complex constants $a_j,b_j$ and $z_j = \exp(ik\ell_j)$. The collection of $k > 0$ for which there is a positive-dimensional $k$-eigenspace is called the {\it spectrum} of $(G,\ell)$. The multiplicity of a real number $k > 0$ is the dimension of the corresponding $k$-eigenspace; an element of the spectra is called {\it simple} if the multiplicity is one and {\it degenerate} otherwise. 
	
	Allow $n$ to denote the number of edges on a graph. The {\it total space} $TG$ of a graph $G$ is standard $2n$-dimensional complex vector space. Any eigenspace embeds in $TG$ by associating to each eigenfunction $\phi$ the $2n$-tuple $$(a_\phi, b_\phi) := (a_1,b_1,\.,a_n,b_n)$$ where the pair $(a_j,b_j)$ refers to the constants in the coordinate representation $\phi_j$. For a choice of eigenvalue $k$ and metric $\ell$, the boundary conditions for the eigenfunctions become linear equations in the coordinates of $TG$, and the matrix of coefficients is called the {\it edge scattering matrix} $S_G(k,\ell)$\textemdash its kernel identifies the $k$-eigenspace of $(G,\ell)$. Each entry of the scattering matrix $S_{G}(k,\ell)$ depends on the functions $z_j : = \exp(ik\ell_j)$ linearly. In this sense, the eigenspace is periodic in the product $k\ell\in \R^n$. If we allow $z = (z_1,\.,z_n)$ to vary, then $S_G$ is a matrix whose entries are polynomials in $z$. Let $$\T^n : = \{ z\in \C^n : |z_j| = 1 \}$$ be the $n$-torus. For any fixed point $z\in \T^n$, the {\it $z$-eigenspace}, $T_{z}G$ is the kernel of $S_G$ evaluated at said point. See [B-K13] for a combinatorial construction of $S_G$.
	
	{\it Remark.} We often prefer to work primarily with $T_{z}G$ rather than an arbitrary pair $k, \ell$. The behavior of eigenfunctions is largely determined by the value $z$ and the corresponding element in $TG$. And as this paper hopes to demonstrate, many non-generic properties code into systems of polynomial equations in $z$.  
	
	For a subgraph $H\subset G$, there is a natural projection $\rho_H$ of total spaces, which remembers coordinates associated to edges in $H$. We abuse notation and write $\rho_H(z), \rho_H(\ell)$ to denote similar projections. In general, there is not a corresponding projection of eigenspaces, but we may still discuss its image. This allows us to characterize `overlap' between eigenspaces of $G$ and a subgraph $H$. We say that a point $z\in \T^n$ is {\it supported} on $H$ if $$\rho_H (T_{z}G)\cap T_{\rho_H(z)} H \neq 0.$$ The statement above is best summarized with an example: Fix some $k$-eigenspace of $(G, \ell)$ such that $z = \exp(ik\ell)$ is supported along $H$, then a nonzero subspace of $k$-eigenfunctions restricts to $k$-eigenfunctions on $H$ with metric $\rho_H(\ell)$. Further, we say $z$ is {\it fully supported} on $H$ if $\rho_H$ defines an inclusion of eigenspaces.
	
	When the subgraph $H$ has no edges, the projection $\rho_H$ does not exist. Instead, we define a so-called {\it evaluation map} $$\eval_{v}(z) : T_zG\to \C$$ for each vertex $v$ and point $z\in \T^n$ by the rule $$\eval_v(z)(a_\phi, b_\phi) := \bigg\{ \begin{array}{ll}
		a_j + b_jz_j & \textnormal{if $e_j\prec v$} \\
		a_jz_j + b_j & \textnormal{if $e_j \succ v$}.
	\end{array}$$ Since elements of $T_zG$ satisfy continuity conditions at vertices, this is a well defined linear map. For a collection $H$ of vertices, let $\eval_H(z)$ denote the direct sum of evaluation maps over the elements in $H$.
	
	{\it III. The Secular Manifold $\Sigma_G$.} The secular polynomial $P_G$ of a graph is the determinant $\det S_G$, so the $z$-eigenspace is positive-dimensional if $z$ is a zero of $P_G$. The collection of zeros forms the {\it secular manifold} $\Sigma_G$, an algebraic subset of $\T^n$. The spectrum of any metric graph $(G, \ell)$ is recovered by intersecting the path $$\gamma_\ell : \R \to \T^n, \;\;\; k\mapsto \exp(ik\ell)$$ with $\Sigma_G$. We sort points of the secular manifold by multiplicity $$\Sigma_G^{(m)} : = \{ z\in \T^n : \dim T_zG \geq m \},$$ or by smoothness, $$\Sigma_G^\reg : = \{ z\in \T^n : z \textnormal{ is regular} \}\;\;\;\;\textnormal{v.s.} \;\;\;\; \Sigma_G^\sing : = \{ z\in \T^n : z \textnormal{ is singular} \}.$$ A theorem of Colin de Verdi\`ere clarifies the relationship between these two categorizations with the equality $\Sigma^\sing_G = \Sigma^{(2)}_G$. Below we state this theorem in a more recognizable and useful format. \\
	
	\noindent {\bf \footnotesize THEOREM 2.2 [CdV].} {\it A point $z$ on $\Sigma_G$ is smooth if and only if $T_{z}G$ is one-dimensional.}
	
	{\it Proof.} Colin de Verdi\`ere's original statement is that an element $k > 0$ of the spectrum of $(G, \ell)$ is simple if and only if $z = \exp(ik\ell)$ is a smooth point of $\Sigma_G$. Since the map $\R_+ \times \R_+^n \to \T^n,$ sending $(k,\ell)\mapsto \exp(ik\ell)$, covers the torus, it may be rephrased in terms of $z$-eigenspaces. $\qed$ \\
	
	Some general information is known about the geometry of secular manifolds. A complete factorization of any secular polynomial is known due to the work of Kurasov \& Sarnak [K-S]. The only example with a nonreduced factor is the circle graph\textemdash its secular polynomial is $(z_1 - 1)^2$. Any factor of a secular polynomial is toral, so its zero locus defines a real analytic hypersurface of the torus. Thus, if we exclude the circle graph, the secular manifold is always reduced with hypersurfaces as irreducible components. Whenever $P_G$ factors, the overlap between any two irreducible components represents a self-intersection of $\Sigma_G$, in which case, it is singular in codimension one [Al]. Colin de Verdi\`ere conjectures that in all other cases, the secular manifold is singular in codimension no less than two [CdV]. 
	
	All relevant statements are summarized in Theorem 2.3, below. At times, it will be useful to consider the zero locus of a secular polynomial $P_G$, as a hypersurface on larger tori. In particular, given an enumeration of the edges of $G$ from a fixed set $\{1,\.,n\}$, we define $Z_G$ to be the zero locus of $P_G$. Note that $Z_G$ is distinct from the secular manifold if $G$ has less than $n$ edges.\\

	\noindent {\bf \footnotesize THEOREM 2.3. [CdV], [K-S]} {\it Let $G$ be a tree graph with enumerated edges $e_{j_1},\.,e_{j_m}$ where $\{ j_1,\.,j_m \}\subset \{ 1,\.,n \}$. The zero locus $Z_G$ of $P_G$ in the $n$-torus is a real analytic hypersurface.}
	
	{\it Proof.} The secular manifold $\Sigma_G\subset \T^m$ is always real analytic since in coordinates, $P_G$ is a trigonometric (analytic) function. The fact that it is a hypersurface follows from the dimension calculation in [CdV] and its irreducibility in [K-S]. Now consider the projection $\T^n\to \T^m$, sending $(z_1,\.,z_n)\mapsto (z_{j_1},\.,z_{j_m})$. It is a submersion, so the preimage $Z_G$ of $\Sigma_G$ is also a real analytic hypersurface (irreducibility follows as before). $\qed$\\
	
	{\it IV. Properties of Generic Eigenfunctions.} The importance of the secular manifold\textemdash from a spectral standpoint\textemdash is most evident in the study of generic eigenfunctions. Some common examples of questions include: {\it Does the metric graph have degenerate spectra for a generic choice of metric?} {\it Do eigenfunctions of a metric graph (for generic metric) satisfy additional boundary conditions (e.g. Dirichlet conditions at a Neumann vertex)?} The answer to any like-minded question is procedural. First, we specify an algebraic subset $V$ of $\Sigma_G$ consisting of points $z$ whose eigenspaces have some undesirable property (e.g. high multiplicity). Next, we consider intersections of $V$ with a generic path of the form $\gamma_\ell$ and compute the expected size of $V\cap \gamma_\ell$. In general, Lior Alon's genericity criterion (Theorem 2.4) provides the answer, while also establishing bounds on just how `generic' a certain property is. For instance, it quickly follows that generic eigenfunctions on graphs are simple and do not satisfy additional boundary conditions, presuming the graph is not a circle. 
	
	We state Theorem 2.4 in the greatest possible generality.\\
	
	\noindent {\bf \footnotesize THEOREM 2.4 [Al].} {\it Let $V$ be a codimension $m$ analytic subvariety of $\T^n$. The set of lengths $$\{ \l\in \R_{+}^n : V\cap \gamma_\ell = \emptyset \}$$ is the complement of a (countable) union of subanalytic sets of codimension at least $m - 1$ in $\R^n_{+}$.} 
	
	{\it Proof.} The bound on the codimension appears directly from Alon's proof of Lemma 6.1 in [Al]. $\qed$ \\
	
	{\it Remark 2.5.} From Theorem 2.4, we deduce that, for most graphs $G$, a generic eigenfunction will not restrict to an eigenfunction of a subgraph. Assume $P_G$ is irreducible, and fix some subgraph $H\subset G$. If $H$ contains all the edges of $G$, then eigenfunctions in $G$ and $H$ must satisfy both Neumann and Dirichlet conditions at some vertex. The collection of $z$ for which $T_zG\cap T_zH$ is nonzero is codimension two in $\T^n$, according to [Al], and hence, non-generic. Suppose now that $H$ is a strict subgraph of $G$. The ideal in $\C[z]$ generated by both $P_G$ and $P_H$ has dimension exactly two, so intersection of $Z_H$ and $\Sigma_G$ must be dimensionally transverse\textemdash that is, of codimension at least two in the torus. Applying Theorem 2.4, for generic $\ell$, none of the eigenfunctions of $(G,\ell)$ restrict to define eigenfunctions on a sub-metric graph. \\
	
	There is no immediate application of Theorem 2.4 once we apply constraints to the possible metrics we wish to consider. A generic metric $\ell$ consists of rationally independent coordinates $\ell_j$, which we express by writing $$\dim \ell := \dim_\Q (\ell_1,\.,\ell_n) = n.$$ The requirement that $\ell$ have a rational span of dimension strictly less than $n$ is non-generic in the parameter space $\R_+^n$. Instead, we specify a subspace $\R_+^m\injec \R_+^n$ cut out by linear equations with coefficients in $\Z$, and consider the induced map $\R_+^m\to \T^n.$ The image is a certain subtorus $M$ of dimension $m$ in $\T^n$. Applying Theorem 2.4 with some algebraic set $V\subset \T^n$ requires dimension bounds on the intersection $V\cap M$. But of course, $V\cap M$ varies with the choice of embedding $\R_+^m\injec \R_n^+$. 
	
	Now consider a fixed graph $G$. For metrics $\ell$ with sufficiently large rational span, the path $\gamma_\ell$ will come arbitrarily close to the singular locus $\Sigma_G^\sing$. In which case, the spectrum of $(G,\ell)$ will not be uniformly discrete\textemdash in particular, the so-called {\it mingap}, $$\textnormal{mg}(G, \ell) : = \lim_{K\to \infty}\inf_{\substack{k > k' > K \\ k, k' \in \textnormal{spec}(G,\ell)}} (k - k'),$$ will be zero. The collection of metrics $\ell$ giving a uniformly discrete spectrum is therefore highly non-generic in $\R_+^n$ and seems heavily dependent on the topology of $G$. For instance, a graph $G$ is a tree if and only if there exists some $\ell$ such that $\mg(G, \ell) > 0$ and $\dim \ell = 2$ [K-S]. \\
	
	\noindent {\bf 3. The Secular Manifold of a Tree.} The main result of this section requires a proper notion of {vanishing} for points on the secular manifold. We say a point $z\in \Sigma_G$ {\it vanishes at a vertex} $v\in G$ if the eigenspace $T_{z}G$ satisfies Dirichlet conditions at $v$, or rather if $\eval_v(z) = 0$. Further, $z$ is said to {\it vanish along an edge} $e_j$ if the $z$-eigenspace is in the kernel of the projection $\rho_{e_j}$. The collection of all vertices and edges for which a fixed point $z$ vanishes along is a closed subgraph of $G$, denoted $G(z)$. The complement $$\supp_G(z) := G \setminus G(z)$$ is an open subgraph of $G$, called the {\it support} of $z$. The space $T_zG$ is always fully supported along $\supp_G(z)$. All these notions agree with our usual expectations for eigenfunctions of metric graphs; that is, a generic $k$-eigenfunction of $(G,\ell)$ vanishes at a vertex or along an edge if and only if the $z = \exp(ik\ell)$ does the same. 
	
	For arbitrary graphs, there seems to be an intimate relationship between various degrees of vanishing and multiplicity. Trees represent an extreme, as evident by the following theorem. Through the remainder of the section, all graphs are assumed to be a tree.\\
	
	\noindent {\bf \footnotesize LEMMA 3.1 (smoothness criteria for trees).} {\it If a point $z$ does not vanish at any vertex of a tree $G$, then $z$ is a smooth point of the secular manifold.}
	
	{\it Proof.} See Appendix A. \\
	
	Nonvanishing has already been established as a generic property of points on secular manifold via analytic [B-K17] and algebraic [Al] methods. It is interesting to see, as Lemma 3.1 so clearly states, that this property belongs solely to smooth points. Its converse gives a necessary condition for degeneracy: a singular point $z\in \Sigma_G$ must vanish somewhere along the graph. Yet, our distinction between nonvanishing v.s. vanishing points of $\Sigma_G$ does not draw the line between smooth and singular. There are many smooth points without full support. In fact, the collection of such points forms a locus in $\Sigma_G$ of codimension one [Al, or see Remark 2.1]. 
	
	The support of $z$ and its vanishing set is closely related to its multiplicity, i.e. the dimension of $T_zG$. The precise statement (Proposition 3.2) is purely topological. \\
	
	\noindent {\bf \footnotesize PROPOSITION 3.2.} {\it Let $z$ be a point of the torus. Then $$\dim T_zG = \beta_0 \; \supp_G(z) - \beta_0\; \partial_G(z).$$ where $\partial_G(z)$ denotes the boundary of $\supp_G(z)$ as a subspace of $G$.} 
	
	{\it Proof.} See Appendix A.\\
	
	We now outline a general approach to calculating components of the singular locus. We say an open subgraph $H\subset G$ is {\it type-m} if $$m = \beta_0 H - \beta_0 \partial_G H.$$ The motivation for this definition appears directly from Proposition 3.2; type-$m$ subgraphs are meant to emulate supports. The idea is to capture all points with support $H.$ Define $$Z(H) : = \bigcap_{K\in \pi_0(H)} Z_K$$ where $Z_K$ is the zero locus of $P_K$ in the $n$-torus. Elements of $Z(H)$ define eigenspaces on each component of $H$. Some glue together to form an $m$-dimensional subspace of $T_zG$, as we prove in Theorem 3.3. It turns out that the components of $\Sigma_G^{(m)}$ are always of this flavor.
	
	As it stands, our definition of a type-$m$ subgraph allows for redundancies. Let $\partial_G H$ denote the boundary of a subgraph $H\subset G$ when considered as a subspace of $G$. If said boundary contains an endpoint $v$ of $H$, then any eigenfunction of $G$ supported on $H$ satisfies both Neumann and Dirichlet conditions at $v$. Necessarily, the eigenfunction vanishes along the entire edge connected at $v$, and $H$ could never be the support of a point $z$ in $\Sigma_G$. To prevent such anomalies, we simply require that $H$ has no endpoints contained in $G$. 
	
	With housekeeping in check, we arrive at the main result: \\
	
	\noindent {\bf \footnotesize THEOREM 3.3.} {\it Let $G$ be a tree graph. The irreducible components of $\Sigma_G^{(m)}$ (with $m > 1$) are the irreducible components of $Z(H)$ where $H$ is an type-$m$ subgraph of $G$. Further, each $Z(H)$ is an algebraic subset of codimension $\beta_0 H$ in $\T^n$.} 
	
	{\it Proof.} Every point of $\Sigma_G^{(m)}$ belongs to some $Z(H)$, with $H$ being type-$m$. Any point $z$ of multiplicity $m$ has its support as an type-$m$ subgraph by Proposition 3.4. Since $z$ must be supported along every component of $\supp_G(z)$, it will lie in $Z(\supp_G(z))$, so there is an inclusion $$\Sigma_G^{(m)}\subset \bigcup_{H\textnormal{ type-}m} Z(H).$$Conversely, we start with an type-$m$ subgraph $H$. Any generic point $z$ of $Z(H)$ will have projections $\rho_K(z)$ nonvanishing along every component $K$ of $H$. The $\rho_H(z)$-eigenspace is therefore $\beta_0H$-dimensional. Applying Neumann conditions at every point in $\partial_G H$ produces an $m$-dimensional subspace of eigenfunctions which extend to the entire graph $G$. The point $z$ will therefore have multiplicity at least $m$ in $\Sigma_G$. Thus, the inclusion (show above) is in fact an equality.
	
	Let's compute $\dim Z(H)$. Since the components $K$ of $H$ are disjoint, the variables occurring in the $P_K$ do not overlap, and the zero loci $Z_K$ must intersect dimensionally transversely. Because each $Z_K$ has codimension one in $\T^n$, $Z(H)$ has codimension $\beta_0 H$.$\qed$ \\
	
	\noindent {\bf \footnotesize COROLLARY 3.4.} {\it The singular locus has codimension exactly two in the secular manifold of a tree with more than two edges and a non-removable singularity. That is, Colin de Verdi\`ere's conjecture holds for trees. All other trees (those with removable singularities) have a smooth secular manifold.}
	
	{\it Proof.} We first establish that the singular locus has dimension no more than $n-3$. Let $H$ be a type-$2$ subgraph of the tree. The boundary must be nonempty, as $G$ is connected; thus, $$\beta_0 H >  \beta_0 H - \beta_0 \partial_G H = 2.$$ By Theorem 3.3, the algebraic set $Z(H)$ must have dimension $n - 3$ or less. Since $H$ was arbitrary, $\dim \Sigma_G^\sing \leq n - 3$. The dimension of $\Sigma_G$ is always $n - 1$ according to Theorem 2.3. The corollary now follows if we prove that $\Sigma_G^\sing$ has a component of dimension $n - 3$. Fix some vertex $v$ of degree greater than two. Let $H$ be a subgraph of $G\setminus v$ consisting of exactly three components; it is type-$2$ with $\beta_0 H = 3$. Applying Theorem 3.3 again completes the equality.
	
	All remaining trees have only removable singularities and are homeomorphic to an interval (maybe missing endpoints). If the graph has $n$ edges, then the secular polynomial takes the form $z_1^2\dotsm z_n^2 \pm 1$, and the gradient is always nonvanishing. $\qed$\\

	\noindent {\bf 4. Applications.} A classic theorem of Lenoid Freidlander is that for a generic choice of metric $\ell$, and for any graph $G$ but the circle, $(G, \ell)$ has simple spectrum [Fr]. Consider the question: {\it does $(G,\ell)$ remain simple for a generic choice of rationally dependent $\ell$}? Neither Friedlander's theorem, nor an improved version found in [Al], can make meaningful conclusions about the spectrum in this case. The reason being that the collection of rationally dependent metrics is non-generic in the parameter space $\R_+^n$. However, when $G$ is a tree, the codimension constraint on the singular locus gives enough flexibility to produce a comparable result for  rationally dependent metrics.
	
	Let $\D$ denote the collection of all $\ell$ with $\dim \ell < n$. We begin by clarifying the geometry of $\mathcal{D}$. Any element satisfies a dependence relation of the form $A_1\ell_1 \+ A_n \ell_n = 0$ for some integer coefficients. Allow $\mathcal{F}_{A}$ to denote the collection of metrics satisfying this constraint. Then $\mathcal{D} = \cup_{A\in \Z^n}\mathcal{F}_{A}.$ A subset $S$ is residual, or subanalytic of positive codimension, in $\mathcal{D}$ if and only if each intersection $S\cap \mathcal{F}_A$ has said property. Then for a statement to hold for ``generic $\ell\in \mathcal{D}$,'' it suffices to work with arbitrary $\mathcal{F}_A$. The conclusion of Theorem 4.1 follows from this line of thinking.\\ 
	
	\noindent {\bf \footnotesize THEOREM 4.1.} {\it Fix a tree graph $G$. Then for a generic choice of rationally dependent $\ell$, $(G, \ell)$ has simple spectrum.}
	
	{\it Proof.} Fix some arbitrary $\mathcal{F} = \mathcal{F}_A$. From the preceeding paragraph, it suffices to check that $(G,\ell)$ is simple for generic $\ell\in \mathcal{F}$. Every path $\gamma_\ell$ is in the subtorus $Z^{A}\subset \T^n$ defined by the equation $z^{A_1}\dotsm z^{A_n} = 1$. Then define $$W_G = Z^A\cap \Sigma_G \;\;\;\; W_G^\sing = Z^A\cap \Sigma_G^\sing.$$ The set $W_G$ forms a real analytic hypersurface of $Z^A$, while $W_G^\sing$ is real analytic of codimension at least two in $Z^A$ by Corollary 3.4. Applying Theorem 2.4 indicates that the set $$\{  \ell\in \mathcal{F} : W_G^\sing \cap \gamma_\ell = \emptyset\}$$ is the complement of a countable union of hypersurfaces. Therefore, for a generic choice of $\ell$ in the family, $(G,\ell)$ has simple spectrum.$\qed$\\
	
	{\it Remark 4.2.} Any graph $G$ with a singular locus of codimension two in $\Sigma_G$ will have the same behavior indicated by Theorem 4.1, since the proof requires nothing more than the codimension constraint. In particular, every graph subject to Colin de Verdi\`ere's conjecture should behave in this manner. \\
	
	We now turn our attention away from metrics with large rational span and toward the other extreme: metrics with a rational span of dimension three or less. Rather than examining the behavior of $(G, \ell)$ for generic $\ell$, we are interested in the possible metrics for which $(G,\ell)$ is uniformly discrete. Recall that a graph is {\it uniformly discrete} when $\mg(G,\ell)$ is nonzero, or (equivalently) when the path $\gamma_\ell$ does not have a limit point on $\Sigma_G^\sing$ [K-S, see Section 2]. The inquiry into whether $\mg(G,\ell)$ is zero or not can thus be rephrased as an intersection. The closure of any path $\gamma_\ell$ defines a subtorus with dimension, $\dim \ell$. A description of $\bar {\gamma_\ell}$ is fairly easy to determine; namely, if $\ell$ lies in $\F_A$, then $\gamma_\ell$ lies in the hypertori defined by $z^{A_1}\dotsm z^{A_n} = 1$. Thus, for a fixed tree graph $G$ and metric $\ell$, we have an explicit description of all relevant subvarieties. 
	
	The situation is ideal for taking intersection products in (oriented) cohomology. The ring $H^\ast (\T^n; \Z)$ is a skew-commutative over $\Z$ with generators $\alpha_1,\.,\alpha_n\in H^1(\T^n, \Z)$. The generator $\alpha_j$ is the class of the hypertorus defined by fixing the value of $z_j$. When two submanifolds $M,N$ are of complementary dimension, the {\it intersection product} $M\cdot N$ is the integer coefficient for $\alpha_1\dotsm \alpha_n$ after taking the product of their respective cohomology classes. If $M$ and $N$ are transverse, $M\cdot N$ is the difference between the positively and negatively oriented points in the intersection $M\cap N$. Thus, the absolute value of $M\cdot N$ always gives a lower bound on $|M\cap N|$, and critically,, $M\cap N$ cannot be empty if $M\cdot N$ is nonzero.
	
	The following theorem forms the intersection theory basis for determining if a metric graph is uniformly discrete or not. As an application, we reconfirm a result of [K-S] using purely algebraic methods. \\

	\noindent {\bf \footnotesize THEOREM 4.3.} {\it If $(G,\ell)$ is uniformly discrete, then the intersection products $\bar {\gamma_\ell}\cdot Z$ are zero for each component $Z$ of the singular locus $\Sigma_G^\sing$. }
	
	{\it Proof.} If this subtorus meets $Z$ nontransversely, then an intersection must occur, so we may assume transversality holds. Thus, the product $\bar{\gamma_\ell}\cdot Z$ is the cohomology class of the intersection $\bar{\gamma_\ell}\cap Z = \emptyset$, which must be zero. \\
	
	\noindent {\bf \footnotesize COROLLARY 4.4 [K-S].} {\it Let $G$ be a star graph. If $\dim \ell > 2$, then the metric graph $(G,\ell)$ cannot be uniformly discrete.}
	
	{\it Proof.} The singular locus of a star graph has codimension three in the torus $\T^n$, so it suffices to prove the statement with $\dim \ell = 3$, in which case $\gamma_\ell$ is a 3-dimensional subtori with nonzero cohomology class in $H^3(\T^n;\Z)$. The highest dimensional components of $\Sigma_G^\sing$ appear from type-2 subgraphs of $G$. There are $\binom{n}{3}$ of them, each corresponding to three edges of $G$ with the central vertex missing. Fix an arbitrary type-$2$ subgraph $H$ with edges $e_{j_1},e_{j_2},e_{j_3}$. Then $$Z(H) = \{ z\in \T^n : z_{j_1}^2 \pm 1 = z_{j_2}^2 \pm 1 = z_{j_3}^2 \pm 1 = 0 \}.$$ As evident from the description above, $Z(H)$ is a transverse intersection of subtori. It's cohomology class (up to orientation) is therefore $8\alpha_{j_1}\alpha_{j_2}\alpha_{j_3}$. With $H$ arbitrary, the cohomology classes of the singular locus are a scaling of the basis for $H^{n-3}(\T^n;\Z)$. The products $\bar{\gamma_\ell}\cdot Z(H)$ cannot all be zero, so $\mg(G,\ell) = 0$. $\qed$\\
	
	For general tree graphs, components of the singular locus do not usually exhibit a manifold structure. Some care must be undertaken to define the cohomology classes that Theorem 4.3 presupposes the existence of. According to Proposition 4.6, the zero loci of any secular polynomial is a pseudomanifold (see Definition 4.5) and admits a cohomology class in $H^1(\T^n; \Z)$. Components of $\Sigma_G$ are then well-behaved intersections of pseudomanifolds, and admit cohomology classes.\\
	
	\noindent {\bf \footnotesize DEFINITION 4.5.} A {\it pseudomanifold} as defined in [G-MP], is a compact space $X$ along with a closed subspace $Y$ of codimension two or more such that $X\setminus Y$ is a smooth orientable manifold, dense in $X$. Any pseudomanifold $X$ of a larger ambient smooth manifold $M$ will define an integer cohomology class as follows: the fundamental class of $X$ pushes forward to $M$ where we may take its Poincar\`e dual. \\
	
	\noindent {\bf \footnotesize PROPOSITION 4.6.} {\it Let $G$ be a tree graph with enumerated edges $e_{j_1},\., e_{j_m}$ where $\{ j_1,\.,j_m \}\subset \{1,\.,n\}$. Then the zero locus $Z_G$ of the secular polynomial $P_G$ on the $n$-torus is a  pseudomanifold with cohomology class $$[Z_G] = 2\alpha_{j_1}\+ 2\alpha_{j_m}$$ up to a choice of orientation.}
	
	{\it Proof.} In coordinates on the $n$-torus, $Z_G$ is defined as the zero level set of a trigonometric function. Thus, $Z_G$ is a compact real analytic space; the set of regular points is always dense. Further, $Z_G^\reg$ is the zero level set of a globally defined function on the oriented manifold $\T^n\setminus Z_G^\sing$, and is therefore orientable. The dimension constraint given by Corollary 3.4 completes the proof that $Z_G$ is a pseudomanifold. 
	
	See Appendix for the calculation of $[Z_G]$. $\qed$ \\

	\noindent {\bf \footnotesize EXAMPLE 5.7.} We illustrate the potential applications of this intersection machinery with an explicit example. Let $G$ be the graph shown in Figure 5.8. Say we wish to study metric graphs $(G,\ell)$ with $\dim \ell = 3$ and $\mg(G,\ell) > 0$. The coordinates of any such metric $\ell$ are related by four dependence relations $$\sum A_j \ell_j = \dotsm = \sum D_j \ell_j = 0$$ where the coefficients $A_1,\.,D_7$ are integers. The torus $\bar{\gamma_\ell}$ lies in special position to the singular locus of $\Sigma_G$, so the coefficients $A_1,\.,D_7$ are subject to algebraic constraints. Some\textemdash certainly not all\textemdash can be calculated directly from the main theorems in Sections 3 \& 4.
	
	\begin{center}

		\tikzset{every picture/.style={line width=0.75pt}} %set default line width to 0.75pt        
		
		\begin{tikzpicture}[x=0.75pt,y=0.75pt,yscale=-0.5,xscale=0.5]
			%uncomment if require: \path (0,476); %set diagram left start at 0, and has height of 476
			
			%Straight Lines [id:da3885529878255418] 
			\draw    (145,114) -- (239,183) ;
			\draw [shift={(239,183)}, rotate = 36.28] [color={rgb, 255:red, 0; green, 0; blue, 0 }  ][fill={rgb, 255:red, 0; green, 0; blue, 0 }  ][line width=0.75]      (0, 0) circle [x radius= 3.35, y radius= 3.35]   ;
			\draw [shift={(145,114)}, rotate = 36.28] [color={rgb, 255:red, 0; green, 0; blue, 0 }  ][fill={rgb, 255:red, 0; green, 0; blue, 0 }  ][line width=0.75]      (0, 0) circle [x radius= 3.35, y radius= 3.35]   ;
			%Straight Lines [id:da8091752354618573] 
			\draw    (239,183) -- (358,183) ;
			\draw [shift={(358,183)}, rotate = 0] [color={rgb, 255:red, 0; green, 0; blue, 0 }  ][fill={rgb, 255:red, 0; green, 0; blue, 0 }  ][line width=0.75]      (0, 0) circle [x radius= 3.35, y radius= 3.35]   ;
			%Straight Lines [id:da7620264825740972] 
			\draw    (358,183) -- (452,252) ;
			\draw [shift={(452,252)}, rotate = 36.28] [color={rgb, 255:red, 0; green, 0; blue, 0 }  ][fill={rgb, 255:red, 0; green, 0; blue, 0 }  ][line width=0.75]      (0, 0) circle [x radius= 3.35, y radius= 3.35]   ;
			%Straight Lines [id:da06198925591025306] 
			\draw    (239,183) -- (148,251) ;
			\draw [shift={(148,251)}, rotate = 143.23] [color={rgb, 255:red, 0; green, 0; blue, 0 }  ][fill={rgb, 255:red, 0; green, 0; blue, 0 }  ][line width=0.75]      (0, 0) circle [x radius= 3.35, y radius= 3.35]   ;
			\draw [shift={(239,183)}, rotate = 143.23] [color={rgb, 255:red, 0; green, 0; blue, 0 }  ][fill={rgb, 255:red, 0; green, 0; blue, 0 }  ][line width=0.75]      (0, 0) circle [x radius= 3.35, y radius= 3.35]   ;
			%Straight Lines [id:da28736201532163075] 
			\draw    (449,115) -- (358,183) ;
			\draw [shift={(449,115)}, rotate = 143.23] [color={rgb, 255:red, 0; green, 0; blue, 0 }  ][fill={rgb, 255:red, 0; green, 0; blue, 0 }  ][line width=0.75]      (0, 0) circle [x radius= 3.35, y radius= 3.35]   ;
			%Straight Lines [id:da16112163733528106] 
			\draw    (449,115) -- (543,184) ;
			\draw [shift={(543,184)}, rotate = 36.28] [color={rgb, 255:red, 0; green, 0; blue, 0 }  ][fill={rgb, 255:red, 0; green, 0; blue, 0 }  ][line width=0.75]      (0, 0) circle [x radius= 3.35, y radius= 3.35]   ;
			%Straight Lines [id:da785993285773606] 
			\draw    (540,47) -- (449,115) ;
			\draw [shift={(540,47)}, rotate = 143.23] [color={rgb, 255:red, 0; green, 0; blue, 0 }  ][fill={rgb, 255:red, 0; green, 0; blue, 0 }  ][line width=0.75]      (0, 0) circle [x radius= 3.35, y radius= 3.35]   ;
			
			% Text Node
			\draw (184,120) node [anchor=north west][inner sep=0.75pt]   [align=left] {$\displaystyle e_{1}$};
			% Text Node
			\draw (169,200) node [anchor=north west][inner sep=0.75pt]   [align=left] {$\displaystyle e_{2}$};
			% Text Node
			\draw (288,160) node [anchor=north west][inner sep=0.75pt]   [align=left] {$\displaystyle e_{3}$};
			% Text Node
			\draw (400,191) node [anchor=north west][inner sep=0.75pt]   [align=left] {$\displaystyle e_{4}$};
			% Text Node
			\draw (405,117) node [anchor=north west][inner sep=0.75pt]   [align=left] {$\displaystyle e_{5}$};
			% Text Node
			\draw (507,134) node [anchor=north west][inner sep=0.75pt]   [align=left] {$\displaystyle e_{6}$};
			% Text Node
			\draw (486,54) node [anchor=north west][inner sep=0.75pt]   [align=left] {$\displaystyle e_{7}$};

		\end{tikzpicture}\\
		
		{\bf Figure 5.8.}
	\end{center} 
	
	We begin by computing the highest dimensional components of $\Sigma_G$ and their cohomology classes. According to Theorem 3.4, these correspond to subgraphs
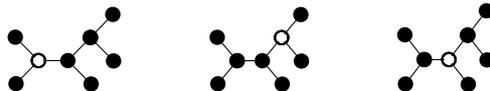
\vspace{-0.5cm} \begin{center}\begin{tikzpicture}[x=0.75pt,y=0.75pt,yscale=-1,xscale=1]
		%uncomment if require: \path (0,476); %set diagram left start at 0, and has height of 476
		
		%Straight Lines [id:da6641574118985116] 
		\draw    (150,339.44) -- (159.93,349.54) ;
		\draw [shift={(161.57,351.22)}, rotate = 45.51] [color={rgb, 255:red, 0; green, 0; blue, 0 }  ][line width=0.75]      (0, 0) circle [x radius= 3.35, y radius= 3.35]   ;
		\draw [shift={(150,339.44)}, rotate = 45.51] [color={rgb, 255:red, 0; green, 0; blue, 0 }  ][fill={rgb, 255:red, 0; green, 0; blue, 0 }  ][line width=0.75]      (0, 0) circle [x radius= 3.35, y radius= 3.35]   ;
		%Straight Lines [id:da3788812774570014] 
		\draw    (163.92,351.22) -- (176.22,351.22) ;
		\draw [shift={(176.22,351.22)}, rotate = 0] [color={rgb, 255:red, 0; green, 0; blue, 0 }  ][fill={rgb, 255:red, 0; green, 0; blue, 0 }  ][line width=0.75]      (0, 0) circle [x radius= 3.35, y radius= 3.35]   ;
		\draw [shift={(161.57,351.22)}, rotate = 0] [color={rgb, 255:red, 0; green, 0; blue, 0 }  ][line width=0.75]      (0, 0) circle [x radius= 3.35, y radius= 3.35]   ;
		%Straight Lines [id:da5432707378273187] 
		\draw    (176.22,351.22) -- (187.8,363) ;
		\draw [shift={(187.8,363)}, rotate = 45.51] [color={rgb, 255:red, 0; green, 0; blue, 0 }  ][fill={rgb, 255:red, 0; green, 0; blue, 0 }  ][line width=0.75]      (0, 0) circle [x radius= 3.35, y radius= 3.35]   ;
		%Straight Lines [id:da48436409983454665] 
		\draw    (159.94,352.91) -- (150.37,362.83) ;
		\draw [shift={(150.37,362.83)}, rotate = 133.98] [color={rgb, 255:red, 0; green, 0; blue, 0 }  ][fill={rgb, 255:red, 0; green, 0; blue, 0 }  ][line width=0.75]      (0, 0) circle [x radius= 3.35, y radius= 3.35]   ;
		\draw [shift={(161.57,351.22)}, rotate = 133.98] [color={rgb, 255:red, 0; green, 0; blue, 0 }  ][line width=0.75]      (0, 0) circle [x radius= 3.35, y radius= 3.35]   ;
		%Straight Lines [id:da049312709108976005] 
		\draw    (187.43,339.61) -- (176.22,351.22) ;
		\draw [shift={(187.43,339.61)}, rotate = 133.98] [color={rgb, 255:red, 0; green, 0; blue, 0 }  ][fill={rgb, 255:red, 0; green, 0; blue, 0 }  ][line width=0.75]      (0, 0) circle [x radius= 3.35, y radius= 3.35]   ;
		%Straight Lines [id:da04460149819297188] 
		\draw    (189.07,341.29) -- (199,351.39) ;
		\draw [shift={(199,351.39)}, rotate = 45.51] [color={rgb, 255:red, 0; green, 0; blue, 0 }  ][fill={rgb, 255:red, 0; green, 0; blue, 0 }  ][line width=0.75]      (0, 0) circle [x radius= 3.35, y radius= 3.35]   ;
		\draw [shift={(187.43,339.61)}, rotate = 45.51] [color={rgb, 255:red, 0; green, 0; blue, 0 }  ][line width=0.75]      (0, 0) circle [x radius= 3.35, y radius= 3.35]   ;
		%Straight Lines [id:da6761983236930835] 
		\draw    (198.63,328) -- (189.06,337.92) ;
		\draw [shift={(187.43,339.61)}, rotate = 133.98] [color={rgb, 255:red, 0; green, 0; blue, 0 }  ][line width=0.75]      (0, 0) circle [x radius= 3.35, y radius= 3.35]   ;
		\draw [shift={(198.63,328)}, rotate = 133.98] [color={rgb, 255:red, 0; green, 0; blue, 0 }  ][fill={rgb, 255:red, 0; green, 0; blue, 0 }  ][line width=0.75]      (0, 0) circle [x radius= 3.35, y radius= 3.35]   ;

	\end{tikzpicture}\hspace{1cm} \begin{tikzpicture}[x=0.75pt,y=0.75pt,yscale=-1,xscale=1]
		%uncomment if require: \path (0,476); %set diagram left start at 0, and has height of 476
		
		%Straight Lines [id:da3885529878255418] 
		\draw    (125,204.44) -- (134.92,216.22) ;
		\draw [shift={(134.92,216.22)}, rotate = 49.9] [color={rgb, 255:red, 0; green, 0; blue, 0 }  ][fill={rgb, 255:red, 0; green, 0; blue, 0 }  ][line width=0.75]      (0, 0) circle [x radius= 3.35, y radius= 3.35]   ;
		\draw [shift={(125,204.44)}, rotate = 49.9] [color={rgb, 255:red, 0; green, 0; blue, 0 }  ][fill={rgb, 255:red, 0; green, 0; blue, 0 }  ][line width=0.75]      (0, 0) circle [x radius= 3.35, y radius= 3.35]   ;
		%Straight Lines [id:da8091752354618573] 
		\draw    (134.92,216.22) -- (147.48,216.22) ;
		\draw [shift={(147.48,216.22)}, rotate = 0] [color={rgb, 255:red, 0; green, 0; blue, 0 }  ][fill={rgb, 255:red, 0; green, 0; blue, 0 }  ][line width=0.75]      (0, 0) circle [x radius= 3.35, y radius= 3.35]   ;
		%Straight Lines [id:da7620264825740972] 
		\draw    (147.48,216.22) -- (157.4,228) ;
		\draw [shift={(157.4,228)}, rotate = 49.9] [color={rgb, 255:red, 0; green, 0; blue, 0 }  ][fill={rgb, 255:red, 0; green, 0; blue, 0 }  ][line width=0.75]      (0, 0) circle [x radius= 3.35, y radius= 3.35]   ;
		%Straight Lines [id:da06198925591025306] 
		\draw    (134.92,216.22) -- (125.32,227.83) ;
		\draw [shift={(125.32,227.83)}, rotate = 129.6] [color={rgb, 255:red, 0; green, 0; blue, 0 }  ][fill={rgb, 255:red, 0; green, 0; blue, 0 }  ][line width=0.75]      (0, 0) circle [x radius= 3.35, y radius= 3.35]   ;
		\draw [shift={(134.92,216.22)}, rotate = 129.6] [color={rgb, 255:red, 0; green, 0; blue, 0 }  ][fill={rgb, 255:red, 0; green, 0; blue, 0 }  ][line width=0.75]      (0, 0) circle [x radius= 3.35, y radius= 3.35]   ;
		%Straight Lines [id:da28736201532163075] 
		\draw    (155.58,206.42) -- (147.48,216.22) ;
		\draw [shift={(157.08,204.61)}, rotate = 129.6] [color={rgb, 255:red, 0; green, 0; blue, 0 }  ][line width=0.75]      (0, 0) circle [x radius= 3.35, y radius= 3.35]   ;
		%Straight Lines [id:da16112163733528106] 
		\draw    (158.59,206.41) -- (167,216.39) ;
		\draw [shift={(167,216.39)}, rotate = 49.9] [color={rgb, 255:red, 0; green, 0; blue, 0 }  ][fill={rgb, 255:red, 0; green, 0; blue, 0 }  ][line width=0.75]      (0, 0) circle [x radius= 3.35, y radius= 3.35]   ;
		\draw [shift={(157.08,204.61)}, rotate = 49.9] [color={rgb, 255:red, 0; green, 0; blue, 0 }  ][line width=0.75]      (0, 0) circle [x radius= 3.35, y radius= 3.35]   ;
		%Straight Lines [id:da785993285773606] 
		\draw    (166.68,193) -- (158.58,202.8) ;
		\draw [shift={(157.08,204.61)}, rotate = 129.6] [color={rgb, 255:red, 0; green, 0; blue, 0 }  ][line width=0.75]      (0, 0) circle [x radius= 3.35, y radius= 3.35]   ;
		\draw [shift={(166.68,193)}, rotate = 129.6] [color={rgb, 255:red, 0; green, 0; blue, 0 }  ][fill={rgb, 255:red, 0; green, 0; blue, 0 }  ][line width=0.75]      (0, 0) circle [x radius= 3.35, y radius= 3.35]   ;	\end{tikzpicture}\hspace{1cm} \begin{tikzpicture}[x=0.75pt,y=0.75pt,yscale=-1,xscale=1]
		%uncomment if require: \path (0,476); %set diagram left start at 0, and has height of 476
		
		%Straight Lines [id:da5773301673116606] 
		\draw    (195,338.09) -- (204.68,350.55) ;
		\draw [shift={(204.68,350.55)}, rotate = 52.13] [color={rgb, 255:red, 0; green, 0; blue, 0 }  ][fill={rgb, 255:red, 0; green, 0; blue, 0 }  ][line width=0.75]      (0, 0) circle [x radius= 3.35, y radius= 3.35]   ;
		\draw [shift={(195,338.09)}, rotate = 52.13] [color={rgb, 255:red, 0; green, 0; blue, 0 }  ][fill={rgb, 255:red, 0; green, 0; blue, 0 }  ][line width=0.75]      (0, 0) circle [x radius= 3.35, y radius= 3.35]   ;
		%Straight Lines [id:da9343641918701124] 
		\draw    (204.68,350.55) -- (214.59,350.55) ;
		\draw [shift={(216.94,350.55)}, rotate = 0] [color={rgb, 255:red, 0; green, 0; blue, 0 }  ][line width=0.75]      (0, 0) circle [x radius= 3.35, y radius= 3.35]   ;
		%Straight Lines [id:da873926596366887] 
		\draw    (218.38,352.4) -- (226.63,363) ;
		\draw [shift={(226.63,363)}, rotate = 52.13] [color={rgb, 255:red, 0; green, 0; blue, 0 }  ][fill={rgb, 255:red, 0; green, 0; blue, 0 }  ][line width=0.75]      (0, 0) circle [x radius= 3.35, y radius= 3.35]   ;
		\draw [shift={(216.94,350.55)}, rotate = 52.13] [color={rgb, 255:red, 0; green, 0; blue, 0 }  ][line width=0.75]      (0, 0) circle [x radius= 3.35, y radius= 3.35]   ;
		%Straight Lines [id:da849214305136863] 
		\draw    (204.68,350.55) -- (195.31,362.82) ;
		\draw [shift={(195.31,362.82)}, rotate = 127.37] [color={rgb, 255:red, 0; green, 0; blue, 0 }  ][fill={rgb, 255:red, 0; green, 0; blue, 0 }  ][line width=0.75]      (0, 0) circle [x radius= 3.35, y radius= 3.35]   ;
		\draw [shift={(204.68,350.55)}, rotate = 127.37] [color={rgb, 255:red, 0; green, 0; blue, 0 }  ][fill={rgb, 255:red, 0; green, 0; blue, 0 }  ][line width=0.75]      (0, 0) circle [x radius= 3.35, y radius= 3.35]   ;
		%Straight Lines [id:da2524391739869234] 
		\draw    (226.32,338.27) -- (218.37,348.68) ;
		\draw [shift={(216.94,350.55)}, rotate = 127.37] [color={rgb, 255:red, 0; green, 0; blue, 0 }  ][line width=0.75]      (0, 0) circle [x radius= 3.35, y radius= 3.35]   ;
		\draw [shift={(226.32,338.27)}, rotate = 127.37] [color={rgb, 255:red, 0; green, 0; blue, 0 }  ][fill={rgb, 255:red, 0; green, 0; blue, 0 }  ][line width=0.75]      (0, 0) circle [x radius= 3.35, y radius= 3.35]   ;
		%Straight Lines [id:da6384872470706267] 
		\draw    (227.76,340.13) -- (236,350.73) ;
		\draw [shift={(236,350.73)}, rotate = 52.13] [color={rgb, 255:red, 0; green, 0; blue, 0 }  ][fill={rgb, 255:red, 0; green, 0; blue, 0 }  ][line width=0.75]      (0, 0) circle [x radius= 3.35, y radius= 3.35]   ;
		\draw [shift={(226.32,338.27)}, rotate = 52.13] [color={rgb, 255:red, 0; green, 0; blue, 0 }  ][line width=0.75]      (0, 0) circle [x radius= 3.35, y radius= 3.35]   ;
		%Straight Lines [id:da7255434218413623] 
		\draw    (235.69,326) -- (227.74,336.41) ;
		\draw [shift={(226.32,338.27)}, rotate = 127.37] [color={rgb, 255:red, 0; green, 0; blue, 0 }  ][line width=0.75]      (0, 0) circle [x radius= 3.35, y radius= 3.35]   ;
		\draw [shift={(235.69,326)}, rotate = 127.37] [color={rgb, 255:red, 0; green, 0; blue, 0 }  ][fill={rgb, 255:red, 0; green, 0; blue, 0 }  ][line width=0.75]      (0, 0) circle [x radius= 3.35, y radius= 3.35]   ;\end{tikzpicture} \end{center}\vspace{-0.4cm} The variety associated to one of these subgraphs is an intersection of zero loci appearing from its connected components. Take for instance, the first subgraph shown above; we have \begin{center}

		\tikzset{every picture/.style={line width=0.75pt}} %set default line width to 0.75pt        
		
		\begin{tikzpicture}[x=0.75pt,y=0.75pt,yscale=-1,xscale=1]
			%uncomment if require: \path (0,476); %set diagram left start at 0, and has height of 476
			
			%Straight Lines [id:da2573612682884123] 
			\draw    (225.22,205.77) -- (232,214.05) ;
			\draw [shift={(232,214.05)}, rotate = 50.69] [color={rgb, 255:red, 0; green, 0; blue, 0 }  ][fill={rgb, 255:red, 0; green, 0; blue, 0 }  ][line width=0.75]      (0, 0) circle [x radius= 3.35, y radius= 3.35]   ;
			\draw [shift={(223.73,203.95)}, rotate = 50.69] [color={rgb, 255:red, 0; green, 0; blue, 0 }  ][line width=0.75]      (0, 0) circle [x radius= 3.35, y radius= 3.35]   ;
			%Straight Lines [id:da25678653535793994] 
			\draw    (231.74,194) -- (225.21,202.12) ;
			\draw [shift={(223.73,203.95)}, rotate = 128.81] [color={rgb, 255:red, 0; green, 0; blue, 0 }  ][line width=0.75]      (0, 0) circle [x radius= 3.35, y radius= 3.35]   ;
			\draw [shift={(231.74,194)}, rotate = 128.81] [color={rgb, 255:red, 0; green, 0; blue, 0 }  ][fill={rgb, 255:red, 0; green, 0; blue, 0 }  ][line width=0.75]      (0, 0) circle [x radius= 3.35, y radius= 3.35]   ;
			%Straight Lines [id:da7358312776643559] 
			\draw    (207.62,213.9) -- (213.38,213.9) ;
			\draw [shift={(215.73,213.9)}, rotate = 0] [color={rgb, 255:red, 0; green, 0; blue, 0 }  ][line width=0.75]      (0, 0) circle [x radius= 3.35, y radius= 3.35]   ;
			\draw [shift={(205.27,213.9)}, rotate = 0] [color={rgb, 255:red, 0; green, 0; blue, 0 }  ][line width=0.75]      (0, 0) circle [x radius= 3.35, y radius= 3.35]   ;
			%Straight Lines [id:da7008513577915432] 
			\draw    (217.22,215.72) -- (224,224) ;
			\draw [shift={(224,224)}, rotate = 50.69] [color={rgb, 255:red, 0; green, 0; blue, 0 }  ][fill={rgb, 255:red, 0; green, 0; blue, 0 }  ][line width=0.75]      (0, 0) circle [x radius= 3.35, y radius= 3.35]   ;
			\draw [shift={(215.73,213.9)}, rotate = 50.69] [color={rgb, 255:red, 0; green, 0; blue, 0 }  ][line width=0.75]      (0, 0) circle [x radius= 3.35, y radius= 3.35]   ;
			%Straight Lines [id:da05872154247101502] 
			\draw    (203.79,215.73) -- (197.26,223.85) ;
			\draw [shift={(197.26,223.85)}, rotate = 128.81] [color={rgb, 255:red, 0; green, 0; blue, 0 }  ][fill={rgb, 255:red, 0; green, 0; blue, 0 }  ][line width=0.75]      (0, 0) circle [x radius= 3.35, y radius= 3.35]   ;
			\draw [shift={(205.27,213.9)}, rotate = 128.81] [color={rgb, 255:red, 0; green, 0; blue, 0 }  ][line width=0.75]      (0, 0) circle [x radius= 3.35, y radius= 3.35]   ;
			%Straight Lines [id:da19923002244193588] 
			\draw    (223.73,203.95) -- (215.73,213.9) ;
			\draw [shift={(215.73,213.9)}, rotate = 128.81] [color={rgb, 255:red, 0; green, 0; blue, 0 }  ][fill={rgb, 255:red, 0; green, 0; blue, 0 }  ][line width=0.75]      (0, 0) circle [x radius= 3.35, y radius= 3.35]   ;
			\draw [shift={(223.73,203.95)}, rotate = 128.81] [color={rgb, 255:red, 0; green, 0; blue, 0 }  ][fill={rgb, 255:red, 0; green, 0; blue, 0 }  ][line width=0.75]      (0, 0) circle [x radius= 3.35, y radius= 3.35]   ;
			%Straight Lines [id:da08582144035441042] 
			\draw    (197,203.8) -- (203.78,212.08) ;
			\draw [shift={(205.27,213.9)}, rotate = 50.69] [color={rgb, 255:red, 0; green, 0; blue, 0 }  ][line width=0.75]      (0, 0) circle [x radius= 3.35, y radius= 3.35]   ;
			\draw [shift={(197,203.8)}, rotate = 50.69] [color={rgb, 255:red, 0; green, 0; blue, 0 }  ][fill={rgb, 255:red, 0; green, 0; blue, 0 }  ][line width=0.75]      (0, 0) circle [x radius= 3.35, y radius= 3.35]   ;
			%Straight Lines [id:da07772414844799314] 
			\draw    (297,216.8) -- (303.78,225.08) ;
			\draw [shift={(305.27,226.9)}, rotate = 50.69] [color={rgb, 255:red, 0; green, 0; blue, 0 }  ][line width=0.75]      (0, 0) circle [x radius= 3.35, y radius= 3.35]   ;
			\draw [shift={(297,216.8)}, rotate = 50.69] [color={rgb, 255:red, 0; green, 0; blue, 0 }  ][fill={rgb, 255:red, 0; green, 0; blue, 0 }  ][line width=0.75]      (0, 0) circle [x radius= 3.35, y radius= 3.35]   ;
			%Straight Lines [id:da6742950793817244] 
			\draw    (361.79,217.73) -- (355.26,225.85) ;
			\draw [shift={(355.26,225.85)}, rotate = 128.81] [color={rgb, 255:red, 0; green, 0; blue, 0 }  ][fill={rgb, 255:red, 0; green, 0; blue, 0 }  ][line width=0.75]      (0, 0) circle [x radius= 3.35, y radius= 3.35]   ;
			\draw [shift={(363.27,215.9)}, rotate = 128.81] [color={rgb, 255:red, 0; green, 0; blue, 0 }  ][line width=0.75]      (0, 0) circle [x radius= 3.35, y radius= 3.35]   ;
			%Straight Lines [id:da615870097370516] 
			\draw    (425.22,219.77) -- (432,228.05) ;
			\draw [shift={(432,228.05)}, rotate = 50.69] [color={rgb, 255:red, 0; green, 0; blue, 0 }  ][fill={rgb, 255:red, 0; green, 0; blue, 0 }  ][line width=0.75]      (0, 0) circle [x radius= 3.35, y radius= 3.35]   ;
			\draw [shift={(423.73,217.95)}, rotate = 50.69] [color={rgb, 255:red, 0; green, 0; blue, 0 }  ][line width=0.75]      (0, 0) circle [x radius= 3.35, y radius= 3.35]   ;
			%Straight Lines [id:da6252103441596972] 
			\draw    (431.74,208) -- (425.21,216.12) ;
			\draw [shift={(423.73,217.95)}, rotate = 128.81] [color={rgb, 255:red, 0; green, 0; blue, 0 }  ][line width=0.75]      (0, 0) circle [x radius= 3.35, y radius= 3.35]   ;
			\draw [shift={(431.74,208)}, rotate = 128.81] [color={rgb, 255:red, 0; green, 0; blue, 0 }  ][fill={rgb, 255:red, 0; green, 0; blue, 0 }  ][line width=0.75]      (0, 0) circle [x radius= 3.35, y radius= 3.35]   ;
			%Straight Lines [id:da5314663302457545] 
			\draw    (407.62,227.9) -- (413.38,227.9) ;
			\draw [shift={(415.73,227.9)}, rotate = 0] [color={rgb, 255:red, 0; green, 0; blue, 0 }  ][line width=0.75]      (0, 0) circle [x radius= 3.35, y radius= 3.35]   ;
			\draw [shift={(405.27,227.9)}, rotate = 0] [color={rgb, 255:red, 0; green, 0; blue, 0 }  ][line width=0.75]      (0, 0) circle [x radius= 3.35, y radius= 3.35]   ;
			%Straight Lines [id:da9998320154682865] 
			\draw    (417.22,229.72) -- (424,238) ;
			\draw [shift={(424,238)}, rotate = 50.69] [color={rgb, 255:red, 0; green, 0; blue, 0 }  ][fill={rgb, 255:red, 0; green, 0; blue, 0 }  ][line width=0.75]      (0, 0) circle [x radius= 3.35, y radius= 3.35]   ;
			\draw [shift={(415.73,227.9)}, rotate = 50.69] [color={rgb, 255:red, 0; green, 0; blue, 0 }  ][line width=0.75]      (0, 0) circle [x radius= 3.35, y radius= 3.35]   ;
			%Straight Lines [id:da9589211394267252] 
			\draw    (423.73,217.95) -- (415.73,227.9) ;
			\draw [shift={(415.73,227.9)}, rotate = 128.81] [color={rgb, 255:red, 0; green, 0; blue, 0 }  ][fill={rgb, 255:red, 0; green, 0; blue, 0 }  ][line width=0.75]      (0, 0) circle [x radius= 3.35, y radius= 3.35]   ;
			\draw [shift={(423.73,217.95)}, rotate = 128.81] [color={rgb, 255:red, 0; green, 0; blue, 0 }  ][fill={rgb, 255:red, 0; green, 0; blue, 0 }  ][line width=0.75]      (0, 0) circle [x radius= 3.35, y radius= 3.35]   ;
			
			% Text Node
			\draw (164,202) node [anchor=north west][inner sep=0.75pt]   [align=left] {$\displaystyle Z$};
			% Text Node
			\draw (176,195) node [anchor=north west][inner sep=0.75pt]  [font=\Large] [align=left] {$\displaystyle ( \ \ \ \ \ \ \ \ )$ \ \ \ \ };
			% Text Node
			\draw (254,201) node [anchor=north west][inner sep=0.75pt]   [align=left] {$\displaystyle =\ Z\ \ \ \ \ \cap \ \ Z\ \ \ \ \ \cap \ \ Z$ };

		\end{tikzpicture}\vspace{-0.4cm}
	\end{center}The strata all intersect each other in the expected dimension, and at regular points, transversely. The cohomology class of the component can thus be written as a product. In particular, the (respective) cohomology classes are $$8\alpha_1\alpha_2(\alpha_3 \+ \alpha_7), \;\;\;\; 8(\alpha_1 \+ \alpha_5)\alpha_6\alpha_7,\;\;\;\; 8(\alpha_1 + \alpha_2 + \alpha_3)\alpha_4(\alpha_5 \+ \alpha_7).$$ For any choice of metric satisfying said dependence relations, the subtorus $\bar{\gamma_\ell}$ is the set of points $$\{ z\in \T^n : z_1^{A_1}\dotsm z_7^{A_7} = \dotsm = z_1^{D_1}\dotsm z_7^{D_7} = 1 \},$$ so its cohomology class is $$(A_1\alpha_1 \+ A_7\alpha_7)\dotsm (D_1\alpha_1 \+ D_7\alpha_7).$$ We now obtain three algebraic expressions in $A_1,\.,D_7$ from taking the product of $[\bar{\gamma_\ell}]$ with the three cohomology classes listed ealier. According to Theorem 4.3, these classes are all zero if $(G,\ell)$ is uniformly discrete. 
	
	These expressions are all too long to include here. Nevertheless, all calculations were preformed with \texttt{Macaulay2}, a computer algebra system. \\

	\noindent {\bf Appendix A.}
	
	\begin{center}

		\tikzset{every picture/.style={line width=0.75pt}} %set default line width to 0.75pt        
		
		\begin{tikzpicture}[x=0.75pt,y=0.75pt,yscale=-.65,xscale=.65]
			%uncomment if require: \path (0,484); %set diagram left start at 0, and has height of 484
			
			%Straight Lines [id:da06708688463532564] 
			\draw    (218,115.5) -- (218,153.5) ;
			\draw [shift={(218,153.5)}, rotate = 90] [color={rgb, 255:red, 0; green, 0; blue, 0 }  ][fill={rgb, 255:red, 0; green, 0; blue, 0 }  ][line width=0.75]      (0, 0) circle [x radius= 3.35, y radius= 3.35]   ;
			\draw [shift={(218,127.5)}, rotate = 90] [color={rgb, 255:red, 0; green, 0; blue, 0 }  ][line width=0.75]    (10.93,-3.29) .. controls (6.95,-1.4) and (3.31,-0.3) .. (0,0) .. controls (3.31,0.3) and (6.95,1.4) .. (10.93,3.29)   ;
			\draw [shift={(218,115.5)}, rotate = 90] [color={rgb, 255:red, 0; green, 0; blue, 0 }  ][fill={rgb, 255:red, 0; green, 0; blue, 0 }  ][line width=0.75]      (0, 0) circle [x radius= 3.35, y radius= 3.35]   ;
			%Straight Lines [id:da688826525187564] 
			\draw    (218,77.5) -- (218,115.5) ;
			\draw [shift={(218,115.5)}, rotate = 90] [color={rgb, 255:red, 0; green, 0; blue, 0 }  ][fill={rgb, 255:red, 0; green, 0; blue, 0 }  ][line width=0.75]      (0, 0) circle [x radius= 3.35, y radius= 3.35]   ;
			\draw [shift={(218,89.5)}, rotate = 90] [color={rgb, 255:red, 0; green, 0; blue, 0 }  ][line width=0.75]    (10.93,-3.29) .. controls (6.95,-1.4) and (3.31,-0.3) .. (0,0) .. controls (3.31,0.3) and (6.95,1.4) .. (10.93,3.29)   ;
			\draw [shift={(218,77.5)}, rotate = 90] [color={rgb, 255:red, 0; green, 0; blue, 0 }  ][fill={rgb, 255:red, 0; green, 0; blue, 0 }  ][line width=0.75]      (0, 0) circle [x radius= 3.35, y radius= 3.35]   ;
			%Straight Lines [id:da934559527225651] 
			\draw    (246,45.5) -- (218,77.5) ;
			\draw [shift={(218,77.5)}, rotate = 131.19] [color={rgb, 255:red, 0; green, 0; blue, 0 }  ][fill={rgb, 255:red, 0; green, 0; blue, 0 }  ][line width=0.75]      (0, 0) circle [x radius= 3.35, y radius= 3.35]   ;
			\draw [shift={(236.61,56.23)}, rotate = 131.19] [color={rgb, 255:red, 0; green, 0; blue, 0 }  ][line width=0.75]    (10.93,-3.29) .. controls (6.95,-1.4) and (3.31,-0.3) .. (0,0) .. controls (3.31,0.3) and (6.95,1.4) .. (10.93,3.29)   ;
			\draw [shift={(246,45.5)}, rotate = 131.19] [color={rgb, 255:red, 0; green, 0; blue, 0 }  ][fill={rgb, 255:red, 0; green, 0; blue, 0 }  ][line width=0.75]      (0, 0) circle [x radius= 3.35, y radius= 3.35]   ;
			%Straight Lines [id:da25121254182778663] 
			\draw    (277,115.5) -- (277,153.5) ;
			\draw [shift={(277,153.5)}, rotate = 90] [color={rgb, 255:red, 0; green, 0; blue, 0 }  ][fill={rgb, 255:red, 0; green, 0; blue, 0 }  ][line width=0.75]      (0, 0) circle [x radius= 3.35, y radius= 3.35]   ;
			\draw [shift={(277,127.5)}, rotate = 90] [color={rgb, 255:red, 0; green, 0; blue, 0 }  ][line width=0.75]    (10.93,-3.29) .. controls (6.95,-1.4) and (3.31,-0.3) .. (0,0) .. controls (3.31,0.3) and (6.95,1.4) .. (10.93,3.29)   ;
			\draw [shift={(277,115.5)}, rotate = 90] [color={rgb, 255:red, 0; green, 0; blue, 0 }  ][fill={rgb, 255:red, 0; green, 0; blue, 0 }  ][line width=0.75]      (0, 0) circle [x radius= 3.35, y radius= 3.35]   ;
			%Straight Lines [id:da0616462898602923] 
			\draw    (277,77.5) -- (277,115.5) ;
			\draw [shift={(277,115.5)}, rotate = 90] [color={rgb, 255:red, 0; green, 0; blue, 0 }  ][fill={rgb, 255:red, 0; green, 0; blue, 0 }  ][line width=0.75]      (0, 0) circle [x radius= 3.35, y radius= 3.35]   ;
			\draw [shift={(277,89.5)}, rotate = 90] [color={rgb, 255:red, 0; green, 0; blue, 0 }  ][line width=0.75]    (10.93,-3.29) .. controls (6.95,-1.4) and (3.31,-0.3) .. (0,0) .. controls (3.31,0.3) and (6.95,1.4) .. (10.93,3.29)   ;
			\draw [shift={(277,77.5)}, rotate = 90] [color={rgb, 255:red, 0; green, 0; blue, 0 }  ][fill={rgb, 255:red, 0; green, 0; blue, 0 }  ][line width=0.75]      (0, 0) circle [x radius= 3.35, y radius= 3.35]   ;
			%Straight Lines [id:da9424756461544965] 
			\draw    (246,45.5) -- (277,77.5) ;
			\draw [shift={(277,77.5)}, rotate = 45.91] [color={rgb, 255:red, 0; green, 0; blue, 0 }  ][fill={rgb, 255:red, 0; green, 0; blue, 0 }  ][line width=0.75]      (0, 0) circle [x radius= 3.35, y radius= 3.35]   ;
			\draw [shift={(256.63,56.47)}, rotate = 45.91] [color={rgb, 255:red, 0; green, 0; blue, 0 }  ][line width=0.75]    (10.93,-3.29) .. controls (6.95,-1.4) and (3.31,-0.3) .. (0,0) .. controls (3.31,0.3) and (6.95,1.4) .. (10.93,3.29)   ;
			\draw [shift={(246,45.5)}, rotate = 45.91] [color={rgb, 255:red, 0; green, 0; blue, 0 }  ][fill={rgb, 255:red, 0; green, 0; blue, 0 }  ][line width=0.75]      (0, 0) circle [x radius= 3.35, y radius= 3.35]   ;
			%Straight Lines [id:da9273884221425932] 
			\draw    (216.59,117.38) -- (196.41,144.12) ;
			\draw [shift={(195,146)}, rotate = 127.02] [color={rgb, 255:red, 0; green, 0; blue, 0 }  ][line width=0.75]      (0, 0) circle [x radius= 3.35, y radius= 3.35]   ;
			\draw [shift={(202.89,135.54)}, rotate = 307.02] [color={rgb, 255:red, 0; green, 0; blue, 0 }  ][line width=0.75]    (10.93,-3.29) .. controls (6.95,-1.4) and (3.31,-0.3) .. (0,0) .. controls (3.31,0.3) and (6.95,1.4) .. (10.93,3.29)   ;
			\draw [shift={(218,115.5)}, rotate = 127.02] [color={rgb, 255:red, 0; green, 0; blue, 0 }  ][line width=0.75]      (0, 0) circle [x radius= 3.35, y radius= 3.35]   ;
			%Straight Lines [id:da42193738094320654] 
			\draw    (277,115.5) -- (308,147.5) ;
			\draw [shift={(308,147.5)}, rotate = 45.91] [color={rgb, 255:red, 0; green, 0; blue, 0 }  ][fill={rgb, 255:red, 0; green, 0; blue, 0 }  ][line width=0.75]      (0, 0) circle [x radius= 3.35, y radius= 3.35]   ;
			\draw [shift={(287.63,126.47)}, rotate = 45.91] [color={rgb, 255:red, 0; green, 0; blue, 0 }  ][line width=0.75]    (10.93,-3.29) .. controls (6.95,-1.4) and (3.31,-0.3) .. (0,0) .. controls (3.31,0.3) and (6.95,1.4) .. (10.93,3.29)   ;
			\draw [shift={(277,115.5)}, rotate = 45.91] [color={rgb, 255:red, 0; green, 0; blue, 0 }  ][fill={rgb, 255:red, 0; green, 0; blue, 0 }  ][line width=0.75]      (0, 0) circle [x radius= 3.35, y radius= 3.35]   ;
			%Straight Lines [id:da07327604638715779] 
			\draw    (308,147.5) -- (302,190) ;
			\draw [shift={(302,190)}, rotate = 98.04] [color={rgb, 255:red, 0; green, 0; blue, 0 }  ][fill={rgb, 255:red, 0; green, 0; blue, 0 }  ][line width=0.75]      (0, 0) circle [x radius= 3.35, y radius= 3.35]   ;
			\draw [shift={(304.16,174.69)}, rotate = 278.04] [color={rgb, 255:red, 0; green, 0; blue, 0 }  ][line width=0.75]    (10.93,-3.29) .. controls (6.95,-1.4) and (3.31,-0.3) .. (0,0) .. controls (3.31,0.3) and (6.95,1.4) .. (10.93,3.29)   ;
			\draw [shift={(308,147.5)}, rotate = 98.04] [color={rgb, 255:red, 0; green, 0; blue, 0 }  ][fill={rgb, 255:red, 0; green, 0; blue, 0 }  ][line width=0.75]      (0, 0) circle [x radius= 3.35, y radius= 3.35]   ;
			%Straight Lines [id:da07981325732338673] 
			\draw    (308,147.5) -- (322,180) ;
			\draw [shift={(322,180)}, rotate = 66.7] [color={rgb, 255:red, 0; green, 0; blue, 0 }  ][fill={rgb, 255:red, 0; green, 0; blue, 0 }  ][line width=0.75]      (0, 0) circle [x radius= 3.35, y radius= 3.35]   ;
			\draw [shift={(312.23,157.32)}, rotate = 66.7] [color={rgb, 255:red, 0; green, 0; blue, 0 }  ][line width=0.75]    (10.93,-3.29) .. controls (6.95,-1.4) and (3.31,-0.3) .. (0,0) .. controls (3.31,0.3) and (6.95,1.4) .. (10.93,3.29)   ;
			\draw [shift={(308,147.5)}, rotate = 66.7] [color={rgb, 255:red, 0; green, 0; blue, 0 }  ][fill={rgb, 255:red, 0; green, 0; blue, 0 }  ][line width=0.75]      (0, 0) circle [x radius= 3.35, y radius= 3.35]   ;
			%Straight Lines [id:da6032232090882677] 
			\draw    (218,153.5) -- (249,185.5) ;
			\draw [shift={(249,185.5)}, rotate = 45.91] [color={rgb, 255:red, 0; green, 0; blue, 0 }  ][fill={rgb, 255:red, 0; green, 0; blue, 0 }  ][line width=0.75]      (0, 0) circle [x radius= 3.35, y radius= 3.35]   ;
			\draw [shift={(228.63,164.47)}, rotate = 45.91] [color={rgb, 255:red, 0; green, 0; blue, 0 }  ][line width=0.75]    (10.93,-3.29) .. controls (6.95,-1.4) and (3.31,-0.3) .. (0,0) .. controls (3.31,0.3) and (6.95,1.4) .. (10.93,3.29)   ;
			\draw [shift={(218,153.5)}, rotate = 45.91] [color={rgb, 255:red, 0; green, 0; blue, 0 }  ][fill={rgb, 255:red, 0; green, 0; blue, 0 }  ][line width=0.75]      (0, 0) circle [x radius= 3.35, y radius= 3.35]   ;
			%Straight Lines [id:da3636796671734297] 
			\draw    (249,185.5) -- (249,223.5) ;
			\draw [shift={(249,223.5)}, rotate = 90] [color={rgb, 255:red, 0; green, 0; blue, 0 }  ][fill={rgb, 255:red, 0; green, 0; blue, 0 }  ][line width=0.75]      (0, 0) circle [x radius= 3.35, y radius= 3.35]   ;
			\draw [shift={(249,210.5)}, rotate = 270] [color={rgb, 255:red, 0; green, 0; blue, 0 }  ][line width=0.75]    (10.93,-3.29) .. controls (6.95,-1.4) and (3.31,-0.3) .. (0,0) .. controls (3.31,0.3) and (6.95,1.4) .. (10.93,3.29)   ;
			\draw [shift={(249,185.5)}, rotate = 90] [color={rgb, 255:red, 0; green, 0; blue, 0 }  ][fill={rgb, 255:red, 0; green, 0; blue, 0 }  ][line width=0.75]      (0, 0) circle [x radius= 3.35, y radius= 3.35]   ;
			%Straight Lines [id:da25785166279133076] 
			\draw    (322,180) -- (336,212.5) ;
			\draw [shift={(336,212.5)}, rotate = 66.7] [color={rgb, 255:red, 0; green, 0; blue, 0 }  ][fill={rgb, 255:red, 0; green, 0; blue, 0 }  ][line width=0.75]      (0, 0) circle [x radius= 3.35, y radius= 3.35]   ;
			\draw [shift={(326.23,189.82)}, rotate = 66.7] [color={rgb, 255:red, 0; green, 0; blue, 0 }  ][line width=0.75]    (10.93,-3.29) .. controls (6.95,-1.4) and (3.31,-0.3) .. (0,0) .. controls (3.31,0.3) and (6.95,1.4) .. (10.93,3.29)   ;
			\draw [shift={(322,180)}, rotate = 66.7] [color={rgb, 255:red, 0; green, 0; blue, 0 }  ][fill={rgb, 255:red, 0; green, 0; blue, 0 }  ][line width=0.75]      (0, 0) circle [x radius= 3.35, y radius= 3.35]   ;
			%Straight Lines [id:da9599176152425319] 
			\draw    (216.61,155.4) -- (196.39,183.1) ;
			\draw [shift={(195,185)}, rotate = 126.14] [color={rgb, 255:red, 0; green, 0; blue, 0 }  ][line width=0.75]      (0, 0) circle [x radius= 3.35, y radius= 3.35]   ;
			\draw [shift={(202.96,174.1)}, rotate = 306.14] [color={rgb, 255:red, 0; green, 0; blue, 0 }  ][line width=0.75]    (10.93,-3.29) .. controls (6.95,-1.4) and (3.31,-0.3) .. (0,0) .. controls (3.31,0.3) and (6.95,1.4) .. (10.93,3.29)   ;
			\draw [shift={(218,153.5)}, rotate = 126.14] [color={rgb, 255:red, 0; green, 0; blue, 0 }  ][line width=0.75]      (0, 0) circle [x radius= 3.35, y radius= 3.35]   ;
			%Straight Lines [id:da3194315633716789] 
			\draw    (248.1,187.67) -- (229.9,231.83) ;
			\draw [shift={(229,234)}, rotate = 112.41] [color={rgb, 255:red, 0; green, 0; blue, 0 }  ][line width=0.75]      (0, 0) circle [x radius= 3.35, y radius= 3.35]   ;
			\draw [shift={(236.71,215.3)}, rotate = 292.41] [color={rgb, 255:red, 0; green, 0; blue, 0 }  ][line width=0.75]    (10.93,-3.29) .. controls (6.95,-1.4) and (3.31,-0.3) .. (0,0) .. controls (3.31,0.3) and (6.95,1.4) .. (10.93,3.29)   ;
			\draw [shift={(249,185.5)}, rotate = 112.41] [color={rgb, 255:red, 0; green, 0; blue, 0 }  ][line width=0.75]      (0, 0) circle [x radius= 3.35, y radius= 3.35]   ;
			%Straight Lines [id:da40091841788660987] 
			\draw    (308,147.5) -- (341,166) ;
			\draw [shift={(341,166)}, rotate = 29.28] [color={rgb, 255:red, 0; green, 0; blue, 0 }  ][fill={rgb, 255:red, 0; green, 0; blue, 0 }  ][line width=0.75]      (0, 0) circle [x radius= 3.35, y radius= 3.35]   ;
			\draw [shift={(329.73,159.68)}, rotate = 209.28] [color={rgb, 255:red, 0; green, 0; blue, 0 }  ][line width=0.75]    (10.93,-3.29) .. controls (6.95,-1.4) and (3.31,-0.3) .. (0,0) .. controls (3.31,0.3) and (6.95,1.4) .. (10.93,3.29)   ;
			\draw [shift={(308,147.5)}, rotate = 29.28] [color={rgb, 255:red, 0; green, 0; blue, 0 }  ][fill={rgb, 255:red, 0; green, 0; blue, 0 }  ][line width=0.75]      (0, 0) circle [x radius= 3.35, y radius= 3.35]   ;
			%Straight Lines [id:da4224511612708055] 
			\draw    (277,153.5) -- (277,191.5) ;
			\draw [shift={(277,191.5)}, rotate = 90] [color={rgb, 255:red, 0; green, 0; blue, 0 }  ][fill={rgb, 255:red, 0; green, 0; blue, 0 }  ][line width=0.75]      (0, 0) circle [x radius= 3.35, y radius= 3.35]   ;
			\draw [shift={(277,165.5)}, rotate = 90] [color={rgb, 255:red, 0; green, 0; blue, 0 }  ][line width=0.75]    (10.93,-3.29) .. controls (6.95,-1.4) and (3.31,-0.3) .. (0,0) .. controls (3.31,0.3) and (6.95,1.4) .. (10.93,3.29)   ;
			\draw [shift={(277,153.5)}, rotate = 90] [color={rgb, 255:red, 0; green, 0; blue, 0 }  ][fill={rgb, 255:red, 0; green, 0; blue, 0 }  ][line width=0.75]      (0, 0) circle [x radius= 3.35, y radius= 3.35]   ;
			%Straight Lines [id:da42461101788671307] 
			\draw    (277,191.5) -- (277,229.5) ;
			\draw [shift={(277,229.5)}, rotate = 90] [color={rgb, 255:red, 0; green, 0; blue, 0 }  ][fill={rgb, 255:red, 0; green, 0; blue, 0 }  ][line width=0.75]      (0, 0) circle [x radius= 3.35, y radius= 3.35]   ;
			\draw [shift={(277,203.5)}, rotate = 90] [color={rgb, 255:red, 0; green, 0; blue, 0 }  ][line width=0.75]    (10.93,-3.29) .. controls (6.95,-1.4) and (3.31,-0.3) .. (0,0) .. controls (3.31,0.3) and (6.95,1.4) .. (10.93,3.29)   ;
			\draw [shift={(277,191.5)}, rotate = 90] [color={rgb, 255:red, 0; green, 0; blue, 0 }  ][fill={rgb, 255:red, 0; green, 0; blue, 0 }  ][line width=0.75]      (0, 0) circle [x radius= 3.35, y radius= 3.35]   ;
			%Straight Lines [id:da5228494883436623] 
			\draw    (277,229.5) -- (299,274) ;
			\draw [shift={(299,274)}, rotate = 63.69] [color={rgb, 255:red, 0; green, 0; blue, 0 }  ][fill={rgb, 255:red, 0; green, 0; blue, 0 }  ][line width=0.75]      (0, 0) circle [x radius= 3.35, y radius= 3.35]   ;
			\draw [shift={(290.66,257.13)}, rotate = 243.69] [color={rgb, 255:red, 0; green, 0; blue, 0 }  ][line width=0.75]    (10.93,-3.29) .. controls (6.95,-1.4) and (3.31,-0.3) .. (0,0) .. controls (3.31,0.3) and (6.95,1.4) .. (10.93,3.29)   ;
			\draw [shift={(277,229.5)}, rotate = 63.69] [color={rgb, 255:red, 0; green, 0; blue, 0 }  ][fill={rgb, 255:red, 0; green, 0; blue, 0 }  ][line width=0.75]      (0, 0) circle [x radius= 3.35, y radius= 3.35]   ;
			%Straight Lines [id:da9523654249989932] 
			\draw    (277,229.5) -- (307,251) ;
			\draw [shift={(307,251)}, rotate = 35.63] [color={rgb, 255:red, 0; green, 0; blue, 0 }  ][fill={rgb, 255:red, 0; green, 0; blue, 0 }  ][line width=0.75]      (0, 0) circle [x radius= 3.35, y radius= 3.35]   ;
			\draw [shift={(296.88,243.75)}, rotate = 215.63] [color={rgb, 255:red, 0; green, 0; blue, 0 }  ][line width=0.75]    (10.93,-3.29) .. controls (6.95,-1.4) and (3.31,-0.3) .. (0,0) .. controls (3.31,0.3) and (6.95,1.4) .. (10.93,3.29)   ;
			\draw [shift={(277,229.5)}, rotate = 35.63] [color={rgb, 255:red, 0; green, 0; blue, 0 }  ][fill={rgb, 255:red, 0; green, 0; blue, 0 }  ][line width=0.75]      (0, 0) circle [x radius= 3.35, y radius= 3.35]   ;
			%Straight Lines [id:da8937115161046638] 
			\draw    (277,264) -- (277,229.5) ;
			\draw [shift={(277,229.5)}, rotate = 270] [color={rgb, 255:red, 0; green, 0; blue, 0 }  ][fill={rgb, 255:red, 0; green, 0; blue, 0 }  ][line width=0.75]      (0, 0) circle [x radius= 3.35, y radius= 3.35]   ;
			\draw [shift={(277,253.75)}, rotate = 270] [color={rgb, 255:red, 0; green, 0; blue, 0 }  ][line width=0.75]    (10.93,-3.29) .. controls (6.95,-1.4) and (3.31,-0.3) .. (0,0) .. controls (3.31,0.3) and (6.95,1.4) .. (10.93,3.29)   ;
			\draw [shift={(277,264)}, rotate = 270] [color={rgb, 255:red, 0; green, 0; blue, 0 }  ][fill={rgb, 255:red, 0; green, 0; blue, 0 }  ][line width=0.75]      (0, 0) circle [x radius= 3.35, y radius= 3.35]   ;
			%Straight Lines [id:da6308282600015145] 
			\draw    (336.07,214.85) -- (336.93,244.65) ;
			\draw [shift={(337,247)}, rotate = 88.34] [color={rgb, 255:red, 0; green, 0; blue, 0 }  ][line width=0.75]      (0, 0) circle [x radius= 3.35, y radius= 3.35]   ;
			\draw [shift={(336.67,235.75)}, rotate = 268.34] [color={rgb, 255:red, 0; green, 0; blue, 0 }  ][line width=0.75]    (10.93,-3.29) .. controls (6.95,-1.4) and (3.31,-0.3) .. (0,0) .. controls (3.31,0.3) and (6.95,1.4) .. (10.93,3.29)   ;
			\draw [shift={(336,212.5)}, rotate = 88.34] [color={rgb, 255:red, 0; green, 0; blue, 0 }  ][line width=0.75]      (0, 0) circle [x radius= 3.35, y radius= 3.35]   ;
			%Straight Lines [id:da8979093753165206] 
			\draw    (336,212.5) -- (358.45,238.23) ;
			\draw [shift={(360,240)}, rotate = 48.89] [color={rgb, 255:red, 0; green, 0; blue, 0 }  ][line width=0.75]      (0, 0) circle [x radius= 3.35, y radius= 3.35]   ;
			\draw [shift={(351.95,230.77)}, rotate = 228.89] [color={rgb, 255:red, 0; green, 0; blue, 0 }  ][line width=0.75]    (10.93,-3.29) .. controls (6.95,-1.4) and (3.31,-0.3) .. (0,0) .. controls (3.31,0.3) and (6.95,1.4) .. (10.93,3.29)   ;
			\draw [shift={(336,212.5)}, rotate = 48.89] [color={rgb, 255:red, 0; green, 0; blue, 0 }  ][fill={rgb, 255:red, 0; green, 0; blue, 0 }  ][line width=0.75]      (0, 0) circle [x radius= 3.35, y radius= 3.35]   ;
			%Straight Lines [id:da06639264653417976] 
			\draw    (247.43,187.25) -- (227.57,209.25) ;
			\draw [shift={(226,211)}, rotate = 132.05] [color={rgb, 255:red, 0; green, 0; blue, 0 }  ][line width=0.75]      (0, 0) circle [x radius= 3.35, y radius= 3.35]   ;
			\draw [shift={(233.48,202.71)}, rotate = 312.05] [color={rgb, 255:red, 0; green, 0; blue, 0 }  ][line width=0.75]    (10.93,-3.29) .. controls (6.95,-1.4) and (3.31,-0.3) .. (0,0) .. controls (3.31,0.3) and (6.95,1.4) .. (10.93,3.29)   ;
			\draw [shift={(249,185.5)}, rotate = 132.05] [color={rgb, 255:red, 0; green, 0; blue, 0 }  ][line width=0.75]      (0, 0) circle [x radius= 3.35, y radius= 3.35]   ;
			%Straight Lines [id:da2815497146787862] 
			\draw    (277,77.5) -- (303,104) ;
			\draw [shift={(303,104)}, rotate = 45.55] [color={rgb, 255:red, 0; green, 0; blue, 0 }  ][fill={rgb, 255:red, 0; green, 0; blue, 0 }  ][line width=0.75]      (0, 0) circle [x radius= 3.35, y radius= 3.35]   ;
			\draw [shift={(277,77.5)}, rotate = 45.55] [color={rgb, 255:red, 0; green, 0; blue, 0 }  ][fill={rgb, 255:red, 0; green, 0; blue, 0 }  ][line width=0.75]      (0, 0) circle [x radius= 3.35, y radius= 3.35]   ;
			%Straight Lines [id:da02126253700862657] 
			\draw    (246,13) -- (246,45.5) ;
			\draw [shift={(246,45.5)}, rotate = 90] [color={rgb, 255:red, 0; green, 0; blue, 0 }  ][fill={rgb, 255:red, 0; green, 0; blue, 0 }  ][line width=0.75]      (0, 0) circle [x radius= 3.35, y radius= 3.35]   ;
			\draw [shift={(246,22.25)}, rotate = 90] [color={rgb, 255:red, 0; green, 0; blue, 0 }  ][line width=0.75]    (10.93,-3.29) .. controls (6.95,-1.4) and (3.31,-0.3) .. (0,0) .. controls (3.31,0.3) and (6.95,1.4) .. (10.93,3.29)   ;
			\draw [shift={(246,13)}, rotate = 90] [color={rgb, 255:red, 0; green, 0; blue, 0 }  ][fill={rgb, 255:red, 0; green, 0; blue, 0 }  ][line width=0.75]      (0, 0) circle [x radius= 3.35, y radius= 3.35]   ;
			
			% Text Node
			\draw (341,194) node [anchor=north west][inner sep=0.75pt]   [align=left] {$\displaystyle v$};
			% Text Node
			\draw (265,42) node [anchor=north west][inner sep=0.75pt]   [align=left] {$\displaystyle e_{j}$};

		\end{tikzpicture}
	\end{center}\begin{adjustwidth}{75pt}{75pt}
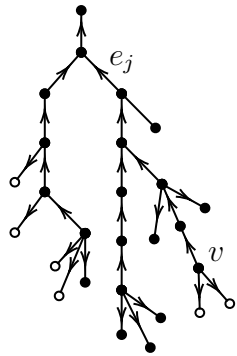
	{\bf Figure A.1.} A tree organized by height, with all edges connecting non-endpoint vertices pointing toward the top. All edges connected to an endpoint are oriented toward them. Here, $v$ is a special vertex below $e_j$.
	\end{adjustwidth}
	
	{\it Proof of Lemma 3.1.} Orient the tree in the manner demonstrated in Figure A.1. A vertex $v$ is said to be {\it special} if all but one edge are connected to endpoints. A vertex $v$ is said to be {\it below} an edge $e_j$ if the vertex is an ancestor of the source of $e_j$. Now, given an element $\phi\in T_zG$, we claim that $\phi_j$ depends linearly on $\eval_v(z)(\phi)$ where $v$ is a special vertex below $e_j$. The proof follows by induction on the maximum distance between $e_j$ and a special vertex below it.
	
	For the base case, $e_j$ has a special vertex $v$ as an endpoint. Set $E = \eval_v(z)(\phi)$. Without loss of generality, we number the edges of $v$ by $e_{n_1},\.,e_{n_p}$ and $e_{d_1},\.,e_{d_q}$ where the first set of edges has Neumann endpoints and the other has a Dirichlet endpoint. Thus, $\phi$ satisfies $$a_{n_m}z_{n_m} = b_{n_m} \;\;\;\; \& \;\;\;\; a_{d_m} z_{d_m} + b_{d_m} = 0,$$ which implies $$E = a_{n_m}(z_{n_m}^2 + 1) \;\;\;\; \& \;\;\;\; E = a_{d_m}(1 - z_{d_m}^2).$$ Since $\eval_v(z)$ is not identically zero (and $z$ is fixed), the constants $z_{n_m}^2 \pm 1$ are nonzero. Hence, $a_{n_m}$ and $a_{d_m}$ are uniquely determined by $E$, and consequently the boundary conditions show that $b_{n_m}$ and $b_{d_m}$ are also dependent on $E$. If $e_j$ is Neumann or Dirichlet, then we are done, otherwise $e_j$ is the one edge which is neither. The Neumann conditions at $v$ are $$a_j + b_jz_j = E \;\;\; \& \;\;\;a_j + \sum \frac{E}{z^2_{n_m} + 1} + \frac{E}{1 - z^2_{d_m}}= b_j z_j + \sum \frac{Ez_{n_m}}{z_{n_m}^2 + 1} - \frac{Ez_{n_m}}{1 - z_{d_m}^2 }.$$ We can then solve for $a_j,b_j$ in terms of $E$ separately since $z_j$ is always nonzero. This completes the base case.
	
	We now proceed to the inductive step. Suppose $e_j$ is distance $N+1$ away from any special vertex below it. Let $v,w$ be any two special vertices below $e_j$. If $s$ is the source of $e_j$, then by induction, there are constants $C_v$ and $C_w$ such that $$C_v\cdot \eval_v(z)(\phi) = \eval_s(z)(\phi) = C_w\cdot \eval_w(z)(\phi)$$ for any element $\phi\in T_zG$. Since $z$ is nonvanishing at $s$, the constants $C_v$ and $C_w$ must be nonzero, and $\eval_v(z)$ and $\eval_w(z)$ differ by a fixed scalar. The pair $v,w$ was arbitrary, so it follows that for any $\phi\in T_zG$, $\phi_m$ depends linearly on $\eval_v(z)(\phi)$ for all edges $e_m$ below $e_j$. The only other edges connected to $s$ are therefore connected to endpoints. An argument identical to the base case allows us to deduce that the value of $\phi$ restricted to these edges are controlled by $\eval_v(z)$. The values of $a_j + b_jz_j$ and $a_j - b_jz_j$ are therefore equal to some fixed scalar multiple of $\eval_v(z)(\phi)$, and $\phi_j$ depends linearly on $\eval_v(z)$. 
	
	The claim in Lemma 3.1 is now clear: the function $\eval_z(v):T_zG\to \C$ is injective for any vertex $v\in G$, so $T_zG$ must be one-dimensional. $\qed$\\
	
	{\it Proof of Proposition 3.2.} For brevity, set $S = \supp_G(z)$ and $B = \partial_G(z)$.
	
	Recall that $z$ is always fully supported on $S$, so the projection $\rho_S$ is an embedding of the $z$-eigenspace. The graph $S$ is potentially disconnected. We note the following relation between its eigenspace and that of its components: $$T_{\rho_S(z)}S = \bigoplus_{H\in \pi_0(S)} T_{\rho_H(z)}H.$$ Each point $\rho_H(z)$ cannot vanish anywhere on $H$ since every vertex in $H$ is contained in the support of $z$. The smoothness criteria ensures that each of the spaces in the direct sum decomposition of $T_{\rho_S(z)}S$ is one-dimensional, so that the total dimension count of $T_{\rho_S(z)}$ is just the zeroth betti number of $S$. Points of $T_zG$ satisfy additional boundary conditions appearing from $B$. We claim that these equations form an independent set, in which case Proposition 3.2 would follow. The proof inducts on the size of $B$. 
	
	For the base case, assume that $B$ consists of exactly one vertex $v$. The connected components $H_j$ of $S$ all share the vertex $v$. Let $e_j$ be the edge of $H_j$ with endpoint $v$, and orient it away from $v$. Since Dirichlet conditions are satisfied at $v$, each point $\phi\in T_zG$ satisfies the Neumann condition at $v$, $\sum a_j = 0$ where the variable $a_j$ corresponds in the usual manner to $\phi_j$. The set of $\phi$ in $T_{\rho_S(z)}S$ satisfying this condition is of codimension one, unless the $a_j$-coefficients of each $\phi$ alll vanish. This cannot be true, otherwise $z$ would have support strictly less than $S$, and $T_zG$ has the appropriate dimension. 
	
	For the inductive step, pick one vertex $v$ in $B$ such that all other points of $B$ lie in the same connected component $G'$ of $G\setminus v$. Enumerate the remaining elements of $\pi_0(S)$ by $H_j$. As before, there is an embedding by $\rho_S$ of $T_zG$ into $$T_{\rho_{G'}(z)}G' \oplus \bigoplus T_{\rho_j(z)}H.$$ The dimension of this sum is $$\beta_0(S) - \beta_0(B) + 1$$ by the induction hypothesis applied to $G'$. Elements of $T_zG$ are subject to one additional constraint from the vertex $v$. Arguing as in the base case yields the result. $\qed$\\

	{\it Proof of Proposition 4.6.} Since $Z_G$ has codimension one, the corresponding cohomology class has the form $$[Z_G] = c_1 \alpha_1 \+ c_m \alpha_m$$ for some integers $c_1,\.,c_m$. Allow $Z$ to denote the subtorus of $\T^n$ defined by relations $z_{j_3} = \dotsm = z_{j_m} = 1$. The space $Z_G\cap Z$ is the intersection of $\T^n$ with the vanishing locus of the ideal $$I = ( P_G, z_{j_3} - 1,\., z_{j_m} - 1).$$ According to theorem \#, the ideal $I$ is equivalent to $$(z_1^2z_2^2 \pm 1, z_{j_3 } - 1,\., z_{j_m} - 1)$$ where $z_1^2z_2^2 \pm 1$ is the secular polynomial of the graph $G'$ obtained by contracting all the edges of $G$ except $e_{j_1},e_{j_2}$. The jacobain of $I$ always has rank $m - 2$, so the intersection $Z_G\cap Z$ is transverse. The manifold $Z_G\cap Z$ has cohomology class $$(2\alpha_{j_1} + 2\alpha_{j_2})\cdot \alpha_{j_3}\dotsm \alpha_{j_m},$$ which is also equal to the product $Z_G\cdot Z$ up to sign. Since $Z$ has class $\alpha_{j_3}\dotsm \alpha_{j_m}$, the coefficients $c_1,c_2$ therefore differ with sign from $2$ and $c_1 = c_2$. Because ${j_1},{j_2}$ are arbitrary, the generalized formula follows. $\qed$\\

	\begin{center}\noindent {\bf Acknowledgements.} \end{center}
	
	\noindent I would like to express the deepest gratitude toward Christopher Judge and his student, Lawford Hatcher for their unwavering support over the last year. I'm also deeply indebted to Lior Alon, whose work inspired much of this research\textemdash our conversations marked significant moments in its evolution. I would also like to thank Gregory Berkolaiko and Peter Sarnak for their communications with this project in its beginning stages; I was heavily influenced by their research on the subject.

	Finally, I thank the department of mathematics at both Indiana University and Harvard University for their support. In particular, I want to recognize Dylan Thurston and Philip Wood for organizing the REUs that made this possible. 
	
	This research was funded by the CMSA and NSF grant \#2051032.
	
	\begin{center}\noindent {\bf References.} \end{center}
	
	\begin{enumerate}
		
		\item[{[Al]}] L. Alon, {\it Generic laplacian eigenfunctions on metric graphs}, J. d'Analyse Math\'ematique, 152 (2024), pp. 729-775. 
		
		\item[{[A-V]}] L. Alon \& C. Vinzant, {\it Gap distributions of Fourier quasicrystals with integer weights via Lee-Yang polynomials}. Rev. Mat. Iberoam. 40 (2024), no. 6, pp. 2203–2250.
		
		\item[{[Ba]}] R. Band, {\it The nodal count $\{0,1,2,3,\.\}$ implies the graph is a tree}, Phil. Trans. R. Soc. A.37220120504 (2014).
		
		\item[{[B-K13]}] G. Berkolaiko \& P. Kuchment, {\it Introduction to Quantum Graphs}, Mathematical Surveys and Monographs, AMS, 186 (2013).
		
		\item[{[B-K17]}] G. Berkolaiko \& W. Liu, {\it Simplicity of eigenvalues and non-vanishing of eigenfunctions of a quantum graph}, J. Math. Anal. Appl., 445 (2017), pp. 803–818.
		
		\item[{[B-W]}] G. Berkolaiko \& B. Winn, {\it Relationship between scattering matrix and spectrum of quantum graphs}, Trans. Amer. Math. Soc., 362 (2010), pp. 6261–6277.
		
		\item[{[CdV]}] Y. Colin de Verdi\`ere, {\it Semi-classical measures on quantum graphs and the Gauß map of the determinant manifold}, Annales Henri Poincar\'e, 16 (2015), pp. 347–364.
		
		\item[{[Fr]}] L. Friedlander, {\it Genericity of simple eigenvalues for a metric graph}, Israel J. Math., 146 (2005), pp. 149–156.
		
		\item[{[G-MP]}] M. Goresky \& R. MacPherson, {\it Intersection Homology Theory}, Topology, 19 (1980), pp. 135-162.
		
		\item[{[K-N]}] P. Kurasov \& M. Nowaczyk, {\it Inverse spectral problem for quantum graphs}, J. Phys. A, 38 (2005), pp. 4901–4915.
		
		\item[{[K-S]}] P. Kurasov \& P. Sarnak, {\it The additive structure of the spectrum of the laplacian on a metric graph}, (in preparation).
		
		\item[{[Ł]}] S. Łojasiewicz, {\it On semi-analytic and subanalytic geometry}, Banach Center Publications 34.1 (1995), pp. 89-104.
	\end{enumerate}

\end{document}